\documentclass[onefignum,onetabnum,final]{siamart190516}

\usepackage{graphicx,epstopdf}
\usepackage{subfig}
\usepackage{enumitem} 
\usepackage{tikz}
\usepackage{titlesec}

\usepackage{breqn} 
\usepackage{fancyhdr}
\usepackage{algorithm}
\usepackage{algorithmic}
\usepackage{lipsum}
\usepackage{amsfonts}

\newtheorem{remark}{Remark}

\usepackage{amssymb}
\usepackage{bbm}
\usepackage[numbers,sort&compress]{natbib}

\newcommand{\mbb}{\mathbb}
\newcommand{\mcal}{\mathcal}
\newcommand{\opn}{\operatorname}
\newcommand{\hsp}{\hspace{0.1cm}}
\newcommand{\hspB}{\hspace{0.3cm}}
\newcommand{\pd}{\partial}
\newcommand{\lan}{\left\langle}
\newcommand{\ran}{\right\rangle}

\DeclareMathOperator*{\argmin}{arg\,min}

\newcommand{\tspace}{D}

\newcommand{\fMN}{f_{M}}
\newcommand{\hMN}{h_{M}}

\newcommand{\fMNNx}{f_{M,N_x}}

\newcommand{\dimXi}{d_{\xi}}

\newcommand{\xiMaxMin}{\xi_{\opn{max/min}}}
\newcommand{\xiMax}{\xi_{\opn{max}}}
\newcommand{\xiMin}{\xi_{\opn{min}}}

\newcommand{\xMax}{x_{\opn{max}}}
\newcommand{\xMin}{x_{\opn{min}}}

\newcommand{\Pcons}{P_{\opn{cons}}}
\newcommand{\Acons}{A_{\opn{cons}}}

\usepackage{lmodern}
\usetikzlibrary{decorations.pathreplacing}
\usetikzlibrary{decorations.shapes}
\usetikzlibrary{calc}

\tikzset{decorate sep/.style 2 args=
{decorate,decoration={shape backgrounds,shape=circle,shape size=#1,shape sep=#2}}}

\headers{Positive L2-minimization Moment Method}{N. Sarna}

\title{A Positive and Stable L2-minimization based moment method for the Boltzmann equation of Gas dynamics\thanks{Submitted to the editors xxxx\funding{N.S is supported by the German Federal Ministry for Economic Affairs and Energy (BMWi) in the joint project "MathEnergy - Mathematical Key Technologies for Evolving Energy Grids", sub-project: Model Order Reduction (Grant number: 0324019B).}}}
\author{Neeraj Sarna\thanks{Corresponding author, Max Planck Institute for Dynamics of Complex Technical Systems, Sandtorstr 1, 39106, Magdeburg, Germany, \email{sarna@mpi-magdeburg.mpg.de}.}}

\begin{document}
\maketitle
\begin{abstract}
We consider the method-of-moments approach to solve the Boltzmann equation of rarefied gas dynamics, which results in the following moment-closure problem. Given a set of moments, find the underlying probability density function. The moment-closure problem has infinitely many solutions and requires an additional optimality criterion to single-out a unique solution. Motivated from a discontinuous Galerkin velocity discretization, we consider an optimality criterion based upon L2-minimization. To ensure a positive solution to the moment-closure problem, we enforce positivity constraints on L2-minimization. This results in a quadratic optimization problem with moments and positivity constraints. We show that a (Courant-Friedrichs-Lewy) CFL-type condition ensures both the feasibility of the optimization problem and the L2-stability of the moment approximation. Numerical experiments showcase the accuracy of our moment method. 
\end{abstract}

\section{Introduction}
Due to modeling assumptions, the Euler and the Navier-Stokes equations become inaccurate as a flow deviates significantly from a thermodynamic equilibrium. This motivates one to consider mathematical models that can approximate flows in all regimes of thermodynamic non-equilibrium. One such model is the Boltzmann equation (BE) that govern the evolution of a probability density function (pdf) $f(x,t,\xi)\in\mbb R^+$ and reads
\begin{gather}
\mcal L(f)=0\hspB\text{where}\hspB \mcal L := \pd_t + \xi\cdot\nabla  - Q. \label{BTE}
\end{gather}
Above, $\xi\in\mbb R^{\dimXi}$ is the molecular velocity with $1\leq d_\xi\leq 3$ being the velocity-dimension, $D:=[0,T]$ is the temporal domain with $T > 0$ being the final time, and $\nabla$ represents a gradient in the spatial domain $\Omega\subseteq \mbb R^d$ with $1\leq d\leq 3$ being the space-dimension. The operator $Q$ is the so-called collision operator and models the inter-particle interaction. Furthermore, the transport operator $\pd_t + \xi\cdot\nabla$ models the free-streaming of the gas molecules. Thus, the BE signifies the fact that the pdf changes due to the free-streaming of the gas molecules and the collisions between them. 

In practical applications, one is not interested in the fine details of a pdf but in the macroscopic quantities like density, velocity, temperature, etc. These quantities can be recovered by taking the velocity-moments of the pdf. This motivates the method-of-moments (MOM) approach, where, rather than directly solving the BE, we solve for a finite number of moments of the pdf. The velocity-moments of the BE provide the governing equation for the moments of $f(x,t,\cdot)$, or the so-called moment equations. However, a finite set of moment equations is not closed---the flux term ($\xi\cdot \nabla f$) results in a moment of degree higher than that included in the moment set. Nevertheless, one can close the moment equations by solving the following moment-closure problem.
\begin{gather}
\text{Moment-closure problem: } \text{given a set of moments, find the underlying pdf.}
\end{gather}

There are infinitely many solutions to the moment-closure problem \cite{MeadMaxEntropy}. To single-out a unique solution, one can introduce an optimality criterion by minimizing a strictly convex functional of the pdf. We choose this functional to be the L2-norm of the pdf. Our choice is motivated by a discontinuous Galerkin (DG) discretization of the velocity domain that we interpret as an L2-minimization problem with moment constraints. Later sections provide further clarification. A major drawback of L2-minimization (and so of a DG-discretization \cite{DGVelocity}) is that it does not penalize the negativity of a solution---for the same reason, a Hermite spectral method is also not positivity preserving \cite{Grad1949}. This is undesirable given that we are approximating a pdf that is positive by definition. There is ample numerical and theoretical evidence supporting the claim that a positive solution to a moment-closure problem better approximates a pdf---see the different works on positivity-preserving moment-methods \cite{CaiSpectral,LucDVM,Levermore,Roman35,posPN,EQMOM,ReviewManuel,Malik}. Furthermore, in theory, a negative solution to a moment-closure problem can result in a negative density and temperature, resulting in a breakdown of the solution algorithm. For this reason, we enforce positivity constraints on our L2-minimization problem. This results in a quadratic optimization problem with moments and positivity constraints. 

For the robustness of the algorithm, the feasibility of the quadratic optimization problem is imperative. We show that a CFL-type condition ensures (i) the feasibility of the optimization problem and (ii) the L2-stability of the moment approximation---we insist that stability is crucial in analyzing the convergence of a moment approximation \cite{SarnaJan2020,MomentsToDVM}. A proof for both these properties hinges on relating our moment approximation to a discrete-velocity-method (DVM). We emphasize that our proof is general in the sense that it is independent of the objective functional being minimized to single-out a unique solution to the moment-closure problem. 

Other than L2-minimization (with positivity constraints) one can consider entropy-minimization. A moment approximation based on entropy minimization has several desirable properties like (i) symmetric hyperbolicity, (ii) presence of an H-theorem, and (iii) positivity of solutions to the moment-closure problem---we refer to \cite{Levermore,RomanEff,GrothMaxEntropy,JunkEntropy,Schneider} and the reference therein for a detailed discussion. Despite the favorable theoretical properties, it is challenging to compute an entropy-minimization based closure. A few reasons for this are as follows. Firstly, to perform entropy-minimization one uses Newton iterations where in every iteration one inverts the Hessian of the objective functional. This Hessian (despite the adaptivity of basis proposed in \cite{AdaptiveGraham}) can become severely ill-conditioned---particularly inside shocks and for large moment sets---leading to a slow (or no) convergence of the Newton solver \cite{Roman35,RomanEff}. Secondly, in every Newton iteration, one needs to compute integrals over the $d_\xi$-dimensional velocity domain. An analytical expression for these integrals is usually unavailable and one seeks a numerical approximation via some quadrature routine. The number of these quadrature points can grow drastically with $d_\xi$, making the solution algorithm expensive for multi-dimensional applications \cite{Roman35,RomanEff,EQMOM}. For instance, the number of tensorized Gauss-Legendre quadrature points grow as $\mcal O(N^{d_\xi})$, where $N$ is the number of quadrature points in one direction.

Replacing entropy minimization by L2-minimization (with positivity constraints) does not necessarily solves the two problems mentioned above.
We use the interior-convex-set algorithm to perform L2-minimization and even for problems with strong shocks and large moment sets, we did not encounter issues with the conditioning of the Hessian. Our results suggest that L2-minimization could be an alternative to entropy-minimization for flow regimes where entropy-minimization losses robustness. Furthermore, since L2-minimization is robust for large moment sets, it is appealing for an adaptive approach where depending upon the accuracy requirements, the moment set can change locally in the space-time domain \cite{MalikAdaptivity,Torrilhon2017}.

We note that although our L2-minimization procedure is robust, to approximate the integrals, we use tensorized Gauss-Legendre quadrature points in the velocity domain, which, we expect, makes L2-minimization expensive. Specialized quadrature points can make both L2 and entropy minimization efficient \cite{EQMOM}. However, these quadrature points do not guarantee the feasibility of the minimization problem---a property crucial for the robustness of the solution algorithm.  To tackle infeasible problems, one can try regularizing the minimization problem by relaxing the moment constraints \cite{RegGraham}. The use of a specialized quadrature with a regularized minimization problem is an interesting direction to pursue and we plan to consider it in the future. 

We acknowledge that our work draws inspiration from the positive $PN$ closure proposed in \cite{posPN} for the radiative transport equation. Indeed, we solve a similar optimization problem as that solved by the positive $PN$ closure. Nevertheless, our work differs from \cite{posPN} in the following ways. Firstly, unlike the linear isotropic collision operator considered in \cite{posPN}, we consider the non-linear Boltzmann-BGK operator that we discretize using entropy-minimization to ensure mass, momentum and energy conservation. Secondly, using the solution of the optimization problem, to close the moment system, authors in \cite{posPN} perform a spherical harmonics based velocity reconstruction of the pdf. Our framework suggests that such a reconstruction is not needed if one uses the same quadrature points to compute moments in the moment equations and to solve the minimization problem. Thirdly, through numerical experiments, we study the convergence of our moment approximation and compare it to the DVM. These studies were not performed in \cite{posPN}. Lastly, we establish robust (under vanishing Knudsen limit) L2-stability estimates for our moment approximation. Let us emphasize that to the best of our knowledge, for gas transport applications, none of the previous works consider L2-minimization based moment-closures with positivity constraints.    

The rest of the article is organized as follows. In \Cref{sec: mom approx} we discuss our moment approximation and the details of the BE. In \Cref{sec: xt disc} we discuss the space-time discretization of the moment equations, the feasibility of the optimization problem, and the stability of the moment approximation. In \Cref{sec: extd multi} we extend our framework to multi-dimensional problems, and in \Cref{sec: num exp} we perform numerical experiments.

\section{Moment Approximation}\label{sec: mom approx}
Throughout this section we consider a one dimensional space-velocity domain i.e., $d=d_\xi=1$ in \eqref{BTE}. An extension to multi-dimensions is straightforward and is discussed in \Cref{sec: extd multi}. We start by discussing a positive L2-minimization based moment-closure and use it later to define a moment approximation for the BE.

\subsection{A positive L2-method-of-moments (pos-L2-MOM)}\label{sec: posL2}
Consider a $m$-th order polynomial in $\xi$ given as $p_m(\xi) := \xi^{m}$. Collect all the different $p_m(\xi)$ upto some order $M\in\mbb N$ in a vector $P_M(\xi)$ given as
\begin{gather}
P_M(\xi):=\left(p_0(\xi),\dots,p_{M-1}(\xi)\right)^T, \label{def PM}
\end{gather}
where $(\cdot)^T$ represents the transpose of a vector. For a function $\xi\mapsto g(\xi)\in\mbb R$, we introduce the shorthand notation
\begin{gather}
\lan g\ran:=\int_{\mbb R}g(\xi)d\xi. 
\end{gather}
Note that the definition of $P_M$ implies that the vector $\lan P_M g \ran$ contains the first $M$ moments of $g$. 

For some moment vector $\lambda \in\mbb R^M$, consider the mathematical formulation of the moment-closure problem described earlier in the introduction
\begin{gather}
\text{Find $g_M$}\hsp :\hsp \lan P_M g_M\ran = \lambda.
\end{gather}
Even for a realizable moment vector $\lambda$ (i.e., there exists a $g^*>0$ such that $\lambda = \lan P_Mg^*\ran $), the above problem can have infinitely many solutions \cite{MeadMaxEntropy}. To single-out a unique solution, we use L2-minimization as an additional optimality criterion. Since L2-minimization does not penalize negativity and since we prefer a positive solution to the moment-closure problem, we explicitly enforce a positivity constraint. This result in an optimization problem given as 
\begin{gather}
g_M:= \argmin_{g^*\in L^2(\mbb R)}\frac{1}{2}\|g^*\|^2_{L^2(\mbb R)}\hsp :\hsp \lan P_M g^*\ran = \lambda,\hsp g^* > 0. \label{L2 min cont}
\end{gather}

In the above minimization problem, as yet, it is unclear how to enforce the positivity constraint almost everywhere on $\mbb R$. To tackle this problem, we consider the following two steps---we refer to~\cite{Roman35,RomanEff,posPN} for similar steps related to the minimum-entropy closure and the positive PN closure.
\begin{enumerate}
\item \textit{Truncate the velocity domain:} We truncate the velocity domain $\mbb R$ to $\Omega_\xi :=[\xiMin,\xiMax]$. A decent estimate for $\xiMaxMin$ follows from the velocity and the temperature field of the gas and is discussed later in \Cref{sec: compute xiMax}. The same sub-section discusses the pros and cons associated with truncating the velocity domain.
\item \textit{Positivity constraints on quadrature points:} To perform the integrals in the minimization problem, we use some quadrature points defined over $\Omega_\xi$.  We enforce the positivity constraints only over these quadrature points. We consider $N$ Gauss-Legendre quadrature points and we denote their weights and abscissas by $\{\omega_i\}_i$ and $\{\xi_i\}_i$, respectively. Using the quadrature points, for some function $\xi\mapsto g(\xi)\in\mbb R$, we define 
\begin{gather}
\lan g\ran \approx \lan g\ran_N := \sum_{i=1}^N\omega_i g(\xi_i).  \label{num approx int}
\end{gather}
For convenience, with $W(g)\in\mbb R^N$ we represent a vector that collects all the values of $g$ at the quadrature points i.e.,
\begin{gather}
\left(W(g)\right)_i := g(\xi_i),\hspB \forall i\in \{1,\dots,N\}.
\end{gather}
\end{enumerate}

With the above two simplifications, the optimization problem in \eqref{L2 min cont} transforms to an optimization problem for $W(g_M)$ given as
\begin{gather}
W(g_M) = \argmin_{W^*\in\mbb R^N}\frac{1}{2}\|W^*\|^2_{l^2}\hsp :\hsp  \underline{A L W^*=\lambda},\hsp W^* > 0. \label{l2 quad}
\end{gather}
To write down the moment constraint (the underlined term) in the above problem, we have used the relation 
\begin{gather}
\lan P_M g^*\ran \approx \lan P_M g^*\ran_N = A L W(g^*), \label{mat-vec form}
\end{gather}
where the matrices $A\in\mbb R^{M\times N}$ and $L\in\mbb R^{N\times N}$ are given as
\begin{gather}
A := \left(P_M(\xi_1),\dots,P_M(\xi_N)\right),\hspB L_{ij}=\begin{cases}
\omega_i,\hsp &i=j\\
0,\hsp &i\neq j
\end{cases}.\label{def A}
\end{gather}
Thus, $L$ is a diagonal matrix containing the quadrature weights $\{\omega_i\}$ at its diagonal, and $A$ is a Vandermonde matrix. Note that in \eqref{l2 quad}, for notational simplicity, we defined $W^* = W(g^*)$. 

\begin{remark}[A DG discretization]\label{remark: DG}
To see the similarity between the pos-L2-MOM and a DG velocity space discretization and understand our motivation behind performing L2-minimization, consider the optimization problem
\begin{gather*}
g^{DG}_M:= \argmin_{g^*\in L^2(\Omega_\xi)}\frac{1}{2}\|g^*\|^2_{L^2(\Omega_\xi)}\hsp :\hsp \int_{\Omega_\xi}P_M g^*d\xi = \lambda. 
\end{gather*}
The above problem is a continuous-in-velocity analogue of \eqref{L2 min cont} but without positivity constraints. Using the first order-optimality conditions, one can conclude that a solution to the above problem is given as (see page-2611 of \cite{posPN} for a similar proof related to the PN closure)
$$g^{DG}_M =\alpha^T P_M, $$
where $\alpha$ is a vector of expansion coefficients related linearly to the moment vector $\lambda$---the exact form of $\alpha$ is not important here. The above expansion is the same as the DG velocity discretization proposed in \cite{DGVelocity}. Thus, one can interpret pos-L2-MOM as a DG velocity discretization with positivity constraints. Note that the DG discretization is not necessarily positive on the quadrature points. Numerical experiments will provide further details.
\end{remark}

\begin{remark}[A Hermite expansion]
One can also interpret a Hermite approximation to a pdf as a solution to a weighted L2-minimization problem---we refer to \cite{Grad1949,Cai2014} for an exhaustive discussion on Hermite expansions. Let $p_m(\xi)$ denote the $m$-th order Hermite polynomial $He_m(\xi)$. Normalize the Hermite polynomials such that they are orthogonal under the inner-product of the weighted L2-space $L^2(\mbb R,\exp(-\xi^2/2))$. Let $P_M$ be as defined in \eqref{def PM}. Note that instead of monomials, the vector $P_M(\xi)$ now contains Hermite polynomials.

Consider a weighted L2-minimization problem given as 
\begin{gather*}
g^{H}_M:= \argmin_{g^*\in L^2(\mbb R,\exp(\xi^2/2))}\frac{1}{2}\|g^*\|^2_{L^2(\mbb R,\exp(\xi^2/2))}\hsp :\hsp \lan P_M g^*\ran = \lambda. 
\end{gather*}
Note that as compared to a DG approximation, in the above optimization problem, we did not truncate the velocity domain. One can show that the solution to the above minimization problem is given as 
$$g_M^H = \lambda^T P_M \exp(-\xi^2/2),$$ 
which is similar to the Hermite spectral method proposed in \cite{Grad1949}. Using the same methodology as for the L2-minimization, one can impose positivity constraints in the above minimization problem and enforce them on a set of Gauss-Hermite quadrature points. We leave the development of a positive weighted L2-minimization based moment method as a part of our future work. 
\end{remark}
\subsubsection{Feasibility of the positive L2-minimization}
It is straightforward to conclude the following. If there exists a $z >0$ such that $\lambda = A Lz$ then the optimization problem in \eqref{l2 quad} is feasible with the feasible point $W^* = z$. We collect this simple, but noteworthy, result as follows. We first define a set of realizable moments 
\begin{gather}
R :=\{\lambda\hsp :\hsp \lambda \in\mbb R^M,\hsp \lambda =AL z,\hsp z> 0\}. \label{def R}
\end{gather}
Using $R$, we collect our statement related to the feasibility of the optimization problem.
\begin{lemma}[Feasibility of the optimization problem] \label{lemma: feasibility}
The optimization problem in \eqref{l2 quad} is feasible if $\lambda\in R$. 
\end{lemma}

Note that for a given $\lambda\in R$, the number of feasible points of the optimization problem vary depending upon the value of $N$ relative to $M$. Let $z>0$ be such that $\lambda = ALz$. A feasible point $W^*$ of the optimization problem \eqref{l2 quad} is a positive solution of the linear system 
$$ALW^* = ALz.$$
Since $AL$ is a full-rank matrix ($A$ is a Vandermonde matrix and the Gauss-Legendre quadrature weights are positive), the above linear system has a unique solution $W^*=z$ for $N\leq M$. Thus, the optimization problem has a single feasible point for $N\leq M$. In contrast, the above linear system has infinitely many positive solutions for $N>M$, resulting in infinitely many feasible points. \footnote{Let $W^*$ be a solution to $ALW^*=ALz$. Let $v$ be an element of the null-space of $AL$---since $AL$ is a flat matrix, its null-space is non-empty. Then, for all $\beta$ such that $\min_i(\beta v_i) > -\min_i(w_i)$, we find that $W^* + \beta v$ is also a feasible point.}

The above discussion indicates that for $N\leq M$, we do not need to perform L2-minimization. A unique positive $W(g_M)$ can be recovered by solving the moment constraint $ALW(g_M)=\lambda$. However, for $N\leq M$, a moment-based approach is meaningless because we can directly compute $W(g_M)$ using a discrete-velocity-method (DVM). Since $N \leq M$, this would be less expensive than first computing $\lambda$ and then computing $W(g_M)$ using the optimization problem. Therefore, in the following discussion we only consider $N > M$. The discussion here becomes clearer when we later relate our moment approximation to a DVM. 

\begin{remark}[Practical considerations while choosing $N$]\label{remark: value N}
Practical considerations suggest a compromise between small and large values of $N$. We use an inter-convex-set algorithm to solve the minimization problem in \eqref{l2 quad}. A crude estimate for the complexity of this algorithm is $\mcal O(N^3)$ \cite{ComplexIntPoint}. Thus, choosing a large value of $N$ increases the computational cost of solving the optimization problem, which, as we discuss later, is the most expensive part of our moment approximation. On the contrary, we do not want $N$ to be so small that the error (measured in some norm) in our moment approximation is dominated by the error in our quadrature approximation. Numerical experiments suggest that choosing $N$ between $2 M$ and $5M$ is a good compromise between accuracy and efficiency.
\end{remark}

\subsection{The Boltzmann Equation (BE)}
Equipping the BE with initial and boundary data provides
\begin{equation}
\begin{gathered}
\mcal L(f) = 0\hsp\text{on}\hsp \Omega\times D\times\mbb R,\hspB f(\cdot,t=0,\cdot)=f_0\hsp\text{on}\hsp\Omega\times\mbb R,\\
 f=f_{in}\hsp\text{on}\hsp \pd\Omega_-\times D. \label{BTE 1D}
\end{gathered}
\end{equation}
Above, the spatial domain is given as $\Omega :=[\xMin,\xMax]$, and $\pd\Omega_-$ is the inflow part of the boundary that reads
\begin{gather}
\pd\Omega_-:=\{(x,\xi)\hsp :\hsp \xi \cdot n(x) \leq 0,\hsp x\in\pd\Omega \}, \label{def OmegaMinus}
\end{gather}
where $n(x)$ is a unit normal at $x\in\pd\Omega$ that points out of the domain. For simplicity, we consider only inflow type boundary conditions and not wall boundary conditions i.e., $f_{in}$ is the given data and is independent of the solution $f$ \cite{Carlos}. An inflow type boundary simplifies our result related to the stability of the moment approximation discussed later. With some additional technical details, one can extend our stability result to solid-wall boundaries---see \cite{Sarna2017} for stability results related to a solid-wall boundary for a Grad's moment method.

We normalise $f$ such that the density $\rho$, the velocity $v$ and the temperature $\theta$ (in energy units) reads
\begin{gather}
\left(\begin{array}{c}
\rho(x,t)\\
\rho(x,t)v(x,t)\\
\rho(x,t)\left(\theta(x,t) + v(x,t)^2\right)
\end{array}\right) := \lan P_{\opn{cons}}f(x,t,\cdot)\ran,\hspB P_{\opn{cons}}(\xi):=\left(\begin{array}{c}
1\\
\xi\\
\xi^2
\end{array}\right). \label{def Pcons}
\end{gather}
Note that for $M\geq 3$, $\Pcons(\xi)$ is nothing but the first three entries of $P_M(\xi)$.

We consider a Boltzmann-BGK collision operator given as 
\begin{gather}
Q(f(x,t,\xi)) := \frac{1}{\tau(x,t)}(f_{\mcal M}(x,t,\xi)-f(x,t,\xi)), \label{BGK op}
\end{gather}
where the collision frequency $\tau(x,t)^{-1}$ reads $\tau(x,t)^{-1}:=C\rho(x,t)\theta(x,t)^{1-\omega}$ with $\omega$ begin the exponent in the viscosity law of the gas \cite{chapman1990}. The collision operator represents the fact that the pdf $f(x,t,\cdot)$ is pushed towards the Maxwell-Boltzmann pdf $f_{\mcal M}(x,t,\cdot)$ given as 
\begin{gather}
f_{\mcal M}(x,t,\xi):= \frac{\rho(x,t)}{\sqrt{2\pi\theta(x,t)}}\exp\left(-\frac{(\xi-v(x,t))^2}{2\theta(x,t)}\right).\label{MB}
\end{gather}
Out of all the pdfs that have the same mass, momentum and energy as $f(x,t,\cdot)$, the pdf $f_{\mcal M}(x,t,\cdot)$ is the one that minimizes the Boltzmann's entropy.  Equivalently,
\begin{gather}
f_{\mcal M}(x,t,\cdot) = \argmin_{f^*(\xi)\geq 0}\left\{\lan f^*\log (f^*)\ran \hsp :\hsp \lan P_{\opn{cons}} f^* \ran =  \lan P_{\opn{cons}}f(x,t,\cdot)\ran \right\}. \label{entropy min}
\end{gather}
Later, we use the above interpretation of $f_{\mcal M}$ to discretize it on a velocity grid. A noteworthy property of $Q(f)$ is its collision invariance i.e., $\lan \Pcons Q(f)\ran = 0$ for all $f$ in the domain of $Q$. This ensures that the BE conserves mass, momentum and energy. By considering $M\geq 3$, which ensures that $\Pcons(\xi)$ is contained in the vector $P_M(\xi)$, and by carefully discretizing the collision operator as in \cite{LucDVM}, we will ensure that our moment system also conserves these quantities.

\subsection{Moment equations}
We present a moment approximation to the BE based upon the pos-L2-MOM described in \Cref{sec: posL2}. To derive a governing equation for the moments $\lan P_M\fMN(x,t,\cdot)\ran_N$, we take (discrete) velocity moments of the BE given in \eqref{BTE 1D} to find
\begin{gather}
\pd_t \lan P_M f(x,t,\cdot)\ran_N + \pd_x\underline{\lan P_M \xi f(x,t,\cdot)\ran_N} = \lan P_M Q(f(x,t,\cdot))\ran_N.
\end{gather}
Recall that $\lan \cdot \ran_N$ is as defined in \eqref{num approx int} and is a numerical approximation to the integral $\lan\cdot \ran$.

The above system of equations is not closed---the underlined flux-term contains a $M$-order moment that is not contained in the moment vector $\lan P_M f(x,t,\cdot)\ran_N$. To close the system of equations, using the moments $\lan P_M f(x,t,\cdot)\ran_N$, we need to approximate the values of $f(x,t,\cdot)$ at the quadrature points i.e., we need to approximate the vector $W(f(x,t,\cdot))$ using the moments $\lan P_M f(x,t,\cdot)\ran_N$. We approximate $W(f(x,t,\cdot))$ by $W(f_M(x,t,\cdot))$. To compute $W(f_M(x,t,\cdot))$, we use the L2-minimization problem given in \eqref{l2 quad} with the moment vector $\lambda$ set to $\lan P_M f(x,t,\cdot)\ran_N$. This results in the following closed set of moment equations
\begin{equation}
\begin{aligned}
\pd_t\lan P_M\fMN \ran_N + \pd_x\lan P_M \xi\fMN \ran_N = &\frac{1}{\tau}(\lan P_M f_{\mcal M,N}\ran_N-\lan P_M f_{M}\ran_N)\hsp\text{on}\hsp \Omega\times D,\\
\hspB \lan P_M\fMN(t=0)\ran_N=&\lan P_Mf_{0}\ran_N\hsp\text{on}\hsp\Omega. \label{mom 1D}
\end{aligned}
\end{equation}
A discretization of the boundary conditions is discussed later during the space-time discretization. Let us emphasis again that to compute the flux term $\lan P_M\xi \fMN(x,t,\cdot)\ran_N$, we only need the value of $W(f_M(x,t,\cdot))$, which is available after solving the L2-minimization problem. The pdf $f_{\mcal M,N}$ is an approximation to $f_{\mcal M}$ and is such that $W(f_{\mcal M,N})$ is a solution to an entropy-minimization problem given as \cite{LucDVM}
\begin{gather}
W(f_{\mcal M,N}(x,t,\cdot)) = \argmin_{W^*\in\mbb R^N_{> 0}}\left\{\sum_i w^*_i\log(w^*_i)\omega_i\hsp :\hsp \Acons L W^* =  \lan\Pcons \fMN(x,t,\cdot)\ran_N \right\}. \label{dis maxwell}
\end{gather}
Note that the moment constraints in the above minimization problem ensure that the moment system \eqref{mom 1D} conserves mass, momentum and energy. Furthermore, the above problem is a discrete-in-velocity analogue of the entropy minimization problem given in \eqref{entropy min}.

\subsection{Computing the velocity cut-off}\label{sec: compute xiMax} Recall that we truncate the velocity domain $\mbb R$ to $\Omega_\xi = [\xiMin,\xiMax]$. We use the same technique as a DVM to compute the velocity cut-off $\xiMaxMin$. The technique is summarised as follows---for further details, we refer to \cite{LocalDVM,LucDVM} and the references therein. Estimating $\xiMaxMin$ using the velocity and the temperature of the gas provides
\begin{equation}
\begin{aligned}
\xiMin:= &\inf_{(x,t)\in\Omega\times D}\left( v(x,t) - c\sqrt{\theta(x,t)}\right),\\
 \xiMax:= &\sup_{(x,t)\in\Omega\times D}\left( v(x,t) + c\sqrt{\theta(x,t)}\right).  \label{def ximax}
\end{aligned}
\end{equation}
From arguments in statistical mechanics, a value of $c$ between $3$ and $4$ is desirable. Choosing $c=3.5$ balances accuracy and computational cost. During numerical experiments, we compare results to a DVM. To ensure that the DVM solution is sufficiently refined, we perform a convergence study by first estimating $\xiMaxMin$ using the initial data and the above formulae and then increasing $\xiMax$ (and decreasing $\xiMin$) till the relative error between two subsequent refinements drops below an acceptable value. We use $\xiMax$ from the last refinement cycle for both the DVM and the pos-L2-MOM---\Cref{sec: num exp} provides further details. In practical applications, for flows that do not show large deviations from thermodynamic equilibrium, one can estimate $v(x,t)$ and $\theta(x,t)$ using a Navier-Stokes solver, which is usually much cheaper than a BE solver \cite{LocalDVM}. 

\begin{remark}[Pros and cons of a space-time-independent $\xiMax$]
Our choice of $\xiMax$ (and $\xiMin$) is space-time-independent, which has both positive and negative consequences. Such a velocity cut-off can be accurate only if, on the entire space-time domain, $f(x,t,\cdot)$ is sufficiently small outside of $\Omega_\xi$. In terms of the macroscopic quantities, we can expect to be accurate only for flows with a velocity and a temperature inside a certain range \cite{LucDVM}. On the positive side, as we discuss later, with a space-time-independent $\xiMax$ it is straightforward to ensure the feasibility of the optimization problem in \eqref{l2 quad}. Furthermore, the stability of the moment equations that we establish later can also be attributed to $\xiMax$ being fixed in space-time. 
\end{remark}

\begin{remark}[A space-time-dependent $\xiMax$]
To overcome the limitations mentioned in the previous remark, similar to \cite{AdaptiveDVM}, one can introduce space-time-dependence in $\xiMax$. We failed to introduce this dependence without sacrificing the feasibility of the optimization problem \eqref{l2 quad} and the stability result discussed later. To overcome the feasibility issue, one can try modifying the optimization problem by regularizing it~\cite{RegGraham}. The regularization adds the moment constraint as a penalty term and tries to minimize both the L2-norm of the pdf and the error in satisfying the moment constraint. As for the stability, it is unclear how one can ensure it with a space-time-dependent $\xiMax$. We leave the development of pos-L2-MOM with space-time adaptive $\xiMax$ as a part of our future work.
\end{remark}

\section{Space-time discretization}\label{sec: xt disc}
\subsection{Preliminaries}
We partition $\Omega=[\xMin,\xMax]$ into $N_x$ intervals given as
\begin{gather}
\Omega=\bigcup_{i=1}^{N_x}\mcal I_i,\hspB \mcal I_i =[x_{i-1/2},x_{i+1/2}],
\end{gather}
where $x_{1/2}=\xMin$ and $x_{N_x+1/2}=\xMax$. With $\{t_i\}_{i=1,\dots,K}\subset\tspace$ we represent a set of discrete time instances such that $0=t_1 < t_2\dots < t_K=T$. For simplicity of notation, we assume that all the space and the time intervals are of the same size $\Delta x$ and $\Delta t$, respectively. An extension to non-uniform space-time grids is straightforward. We denote the finite volume (FV) approximation of $\lan P_M\fMN(x,t,\cdot)\ran$ and $\lan P_Mf_{\mcal M,N}(x,t,\cdot)\ran$ in the $i$-th cell and at the $k$-th time instance by
\begin{equation}
\begin{aligned}
\lan P_Mf_i^k \ran_N  \approx &\frac{1}{\Delta x}\int_{\mcal I_i}\lan P_M\fMN(x,t_k,\cdot)\ran_N dx,\\
 \lan P_Mf_{\mcal M,i}^k \ran_N  \approx &\frac{1}{\Delta x}\int_{\mcal I_i}\lan P_Mf_{\mcal M,N}(x,t_k,\cdot)\ran_N dx.  
\end{aligned}
\end{equation}
Above, $f_{\mcal M,N}$ is the discretization of the Maxwell-Boltzmann distribution introduced in \eqref{entropy min} and for notational simplicity, we have suppressed the $M$ dependence in $f_i^k$. 

Using the matrix $A$ and $L$ given in \eqref{def A}, we can express the space-time discrete moments in a matrix-vector product form as
\begin{equation}
\begin{gathered}
\lan P_M f_i^k\ran_N = A LW(f_i^k),\hspB\lan P_M f_{\mcal M,i}^k\ran = ALW(f_{\mcal M,i}^k),
\end{gathered}
\end{equation}
where $W(f_i^k)$ and $W(f_{\mcal M,i}^k)$ are the FV-approximations to $W(\fMN(x,t_k,\cdot))$ and\\ $W(f_{\mcal M}(x,t_k,\cdot))$, respectively, in the $i$-th cell and at the $k$-th time step. For later convenience, with $f_{M,N_x} $ we represent an FV approximation to $f_{M}$ defined as 
\begin{gather}
f_{M,N_x}(x,t_k,\xi) = f_i^k(\xi),\hspB\forall x\in\mcal I_i,k\in\{1,\dots,K\},\xi\in\{\xi_i\}_i. \label{FV fmn}
\end{gather}

\subsection{Evolution scheme}\label{sec: evo scheme}
The evolution scheme consists of four steps outlined below. We represent these steps for some $t=t_k$. Each step is repeated from $k=1$ to $k = K-1$. For $k=1$, we initialize with 
\begin{gather}
\lan P_M f_i^k\ran_N =  \frac{1}{\Delta x}\int_{\mcal I_i}\lan P_M f_{0}(x,\cdot)\ran_N dx,\hspB\forall i\in \{1,\dots,N_x\}, \label{initialization}
\end{gather}
where $f_{0}$ is the initial data in \eqref{BTE 1D}. We approximate the space integral with 10 Gauss-Legendre quadrature points in each cell.
\begin{enumerate}
\item \textit{Entropy-minimization step: }Using the conserved moments $\{\lan\Pcons f_i^k\ran_N\}_i$, solve the entropy minimization problem in \eqref{dis maxwell}. This provides the discrete Maxwell-Boltzmann pdf $\{f_{\mcal M,i}^k\}_{i}$.
\item \textit{Collision step: }With the output of the previous step, perform collisions with an implicit Euler time-stepping scheme. At an intermediate $t=t_{k^*}$ and for all $i\in\{1,\dots,N_x\}$, this provides \cite{FilbetAP}
\begin{gather}
\frac{\lan P_M f_{i}^{k^*}\ran_N-\lan P_M f_{i}^k\ran_N}{\Delta t} = \frac{1}{\tau(x_i,t_{k^*})}\left(\lan P_M f_{\mcal M,i}^{k^*}\ran_N-\lan P_M f_{i}^{k^*}\ran_N\right). \label{col step}
\end{gather}
There is an explicit solution to the above implicit collision step. Since the collision step preserves mass, moment and energy and since the solution of the entropy-minimization problem \eqref{dis maxwell} is unique for a given set of conserved moments, we find $W(f_{\mcal M,i}^{k^*}) = W(f_{\mcal M,i}^{k})$. This implies that $\lan P_M f_{\mcal M,i}^{k^*}\ran_N = \lan P_M f_{\mcal M,i}^{k}\ran_N$, which provides
\begin{equation}
\begin{aligned}
\lan P_M f_{i}^{k^*}\ran_N = &\frac{1}{1+\Delta t/\tau(x_i,t_{k^*})}\lan P_M f_{i}^{k}\ran_N\\
& + \frac{\Delta t/\tau(x_i,t_{k^*})}{1+\Delta t/\tau(x_i,t_{k^*})}\lan P_M f_{\mcal M,i}^{k}\ran_N.\label{explicit collision}
\end{aligned}
\end{equation}

\item \textit{Optimization step: }Using the moments $\{\lan P_M f_{i}^{k^*}\ran_N\}_{i}$, compute the weights $\{W(f_{i}^{k^*})\}_{i}$ by solving the optimization problem in \eqref{l2 quad}. 
\item \textit{Transport step: }Using the output of the previous step, perform the transport step given as
\begin{equation}
\begin{aligned}
\frac{\lan P_M f_{i}^{k+1}\ran_N-\lan P_M f_{i}^{k^*}\ran_N}{\Delta t} = -\frac{1}{\Delta x}\left(\right.&\left.\mcal F(W(f^{k,*}_{i+1}),W(f^{k,*}_{i}))\right.\\
&-\left.\mcal F(W(f^{k,*}_{i}),W(f^{k,*}_{i-1}))\right). \label{transport}
\end{aligned}
\end{equation}
To impose boundary conditions, for $i=1$, set $W(f^{k,*}_{i-1})=W(f_{in}(t,\cdot))$ and for $i=N_x$, set $W(f^{k,*}_{i+1})=W(f_{in,N}(t,\cdot))$, where $f_{in}$ is the boundary data given in \eqref{BTE 1D}. 
Above, $\mcal F:\mbb R^N\times\mbb R^N\to\mbb R^M$ is the numerical flux and since we consider a kinetic upwind numerical flux, it reads \cite{Malik}
\begin{gather}
\mcal F(W_1,W_2) :=\frac{1}{2}\left(AL(\Xi-|\Xi|)W_1 + AL(\Xi+|\Xi|)W_2\right).\label{num flux}
\end{gather}
Above, $A$ and $L$ are the two matrices defined in \eqref{def A}. The matrix $\Xi$ is a diagonal matrix with the locations of the quadrature points $\{\xi_i\}_i$ at its diagonal. Furthermore, $|\Xi|$ is a matrix representing the absolute value of $\Xi$ in the sense that $(|\Xi|)_{ij} = |\Xi_{ij}|$. For clarity, to express $\mcal F$ in a standard kinetic upwind flux form, note that $AL(\Xi\pm|\Xi|)W_1 = \lan P_M(\xi\pm|\xi|) f_1\ran_N$. 
\end{enumerate}

\begin{remark}[Space-time locality of the optimization step]
The optimization step (and also the entropy minimization step) is a local in space-time operation. We loop over each spatial cell, solve the optimization problem, add the local contributions to the numerical flux and move over to the next cell. Therefore, at any given point in time, we store only the moments in all the spatial cells and not the weights. This results in a drastic reduction in memory consumption since, in practice, the number of weights are much larger than the number of moments---see \cite{Roman35,RomanEff} for a similar comment related to a maximum-entropy closure. Let us emphasis that in comparison, a DVM stores the weights in all the cells, which, particularly for multi-dimensional velocity domain, results in a memory intensive algorithm \cite{LocalDVM}.
\end{remark}

\subsection{Properties of the evolution scheme}
The entropy-minimization problem in \eqref{dis maxwell} ensures that our moment approximation conserves mass, moment and energy. In addition to being conservative, the following discussion establishes that our space-time discrete moment approximation (i) under a CFL-type condition, results in a feasible optimization problem; and (ii) is L2 stable in the sense that the L2-energy $\sum_{i=1}^{N_x}\|\lan P_M f_{i}^k\ran_N\|^2_{l^2}$ has an upper-bound that depends solely on the initial data $f_{0}$ and the boundary data $f_{in}$. 

We start with making the following assumptions on the initial and the boundary data. We assume that the first $M$-moments of $f_0$ and $f_{in}$ belong to the realizability set $R$ defined in \eqref{def R} i.e.,
\begin{equation}
\begin{gathered}
\lan P_M f_{in}(x,t,\cdot)\ran_N \in R,\hspB \lan P_M f_{0}(x,\cdot)\ran_N \in R,\hspB
\forall (x,t)\in\Omega\times D.\label{assume data}
\end{gathered}
\end{equation}
The above assumption will be helpful in establishing the feasibility of the optimization problem in the optimization step. For the boundary data, we also assume that 
\begin{equation}
\begin{gathered}
| f_{in}(\cdot,t,\cdot)|_{\pd\Omega,N}<\infty,\hspB\forall t\in D,\\
\text{where  }|f_{in}(\cdot,t,\cdot)|^2_{\pd\Omega,N} := \sum_{\xi_i \cdot n(x)\leq 0}\oint_{\pd\Omega} |\xi_i \cdot n(x)| f_{in}(x,t,\xi_i)^2\omega_i ds. \label{def bc norm}
\end{gathered}
\end{equation}
Above, the unit vector $n(x)$ is as given in \eqref{def OmegaMinus}, and $\{\xi_i\}$ and $\{\omega_i\}_i$ are the abscissas and the weights of the quadrature points, respectively. Note that the assumption on $|f_{in}(\cdot,t,\cdot)|_{\pd\Omega,N}$ is a discrete-in-velocity analogue of a standard assumption that $f_{in}(\cdot,t,\cdot)\in L^2(\pd\Omega_-,|\xi\cdot n(x)|)$---see \cite{Ukai} for further details. Here, $L^2(\pd\Omega_-,|\xi\cdot n(x)|)$ represents a $L^2$ space over $\pd\Omega_-$ with the Lebesgue measure $|\xi\cdot n(x)|$, and the set $\pd\Omega_-$ contains all the incoming velocities and is as defined in \eqref{def OmegaMinus}. Intuitively, the above assumption states that the total L2-energy flux associated with $f_{in}$ should be bounded. We insist that the assumption is valid for most applications of practical relevance.

\subsection{Feasibility of the optimization problem} We show that under a CFL-condition, the moments resulting from the collision step and the transport step belong to the realizability set $R$ given in \eqref{def R} i.e., both the steps are realizability preserving. The feasibility of the optimization problem then follows from \cref{lemma: feasibility}. The details are as follows.

Our result is a straightforward extension of the proof for the realizability preserving space-time discretization of radiative transport equations considered in \cite{HighOrderGraham}. Using the definition of $R$ given in \eqref{def R}, we find
\begin{gather}
a_1\lambda_1+a_2 \lambda_2\in R,\hspB\forall a_1,a_2\geq 0,\hsp \lambda_1,\lambda_2\in R.\label{pos combo}
\end{gather}
We consider the collision step given in \eqref{explicit collision}. Suppose that 
$
\lan P_M f_i^k\ran\in R,$ which implies that entropy-minimization step is well-posed and that $\lan P_M f_{\mcal M,i}^k\ran\in R.$
Then, the above relation implies that for any $\Delta t,\tau(x_i,t_k) > 0$, the collision step is realizability preserving i.e.,   for all $i\in\{1,\dots,N_x\}$, we have $\lan P_M f_i^{k^*}\ran\in R$.

We show that under a CFL-condition, the transport step in \eqref{transport} is also realizability preserving. Replacing the numerical flux function from \eqref{num flux} in the transport step given in \eqref{transport} and re-arranging a few terms provides
\begin{equation}
\begin{aligned}
\lan P_Mf_i^{k+1}\ran_N = &AL(1-\Lambda|\Xi|)W( f_i^{k^*})\\
& + \underline{\frac{\Lambda}{2}AL(|\Xi|-\Xi)W( f_{i+1}^{k^*})}+\underline{\frac{\Lambda}{2}AL(|\Xi| + \Xi)W(f_{i-1}^{k^*})}.\label{conv step reform}
\end{aligned}
\end{equation}
where $\Lambda := \frac{\Delta t}{\Delta x}$. For all $i\in\{1,\dots,N_x\}$, due to the positivity constraints in the optimization problem \eqref{l2 quad}, we have $W(f_{i}^{k^*})>0$, which, for $\Lambda > 0$, implies that the underlined terms are in $R$. To ensure that the first term on the right is in $R$, we choose 
\begin{gather}
0 < \Lambda \leq \opn{min}\{|\xiMax^{-1}|,|\xiMin^{-1}|\}. \label{constraint lambda}
\end{gather}
The above range of $\Lambda$, the relation in \eqref{pos combo} and the assumption on the initial and the boundary data \eqref{assume data} provides $\lan P_Mf_i^{k+1}\ran_N\in R$. We collect our findings in the result below.
\begin{lemma}
Consider the evolution scheme outlined in \Cref{sec: evo scheme} and define $\Lambda = \Delta t/\Delta x$. Assume that the initial and the boundary data satisfies \eqref{assume data}, then the quadratic optimization problem in the evolution scheme is feasible if $\Lambda\in (0,\opn{min}\{|\xiMax^{-1}|,|\xiMin^{-1}|\}]$.
\end{lemma}
\subsection{$L2$ stability of the scheme}
Define the total L2-energy at $t=t_{k+1}$ as
\begin{gather}
\mcal E_{k+1} := \sum_{i=1}^{N_x}\|\lan P_Mf_i^{k+1}\ran_N\|_{l^2}^2. 
\end{gather}
We establish that $\mcal E_{k+1}$ is bounded by the L2-energy of the previous time-step $\mcal E_{k}$ and $|f_{in}(\cdot,t_k,\cdot)|_{\pd\Omega,N}$. Recursion then implies that $\mcal E_{k+1}$ is bounded solely by the initial and the boundary data. 

For convenience, we define a few objects. For a vector $z\in\mbb R^N$, with $\|z\|_L$ we represent the norm
$$
\|z\|_L := \sqrt{z^T L z}. 
$$
Interpreting $z$ as a vector that contains the value of a function $g:\mbb R\to\mbb R$ at the quadrature points, we conclude that $\|z\|_L$ represent an approximation to $\|g\|_{L^2}$. We bound the $l^2$-norm of a moment vector $\lambda = ALz$ as 
\begin{gather}
\|\lambda\|_{l^2}\geq \sigma_\min (A\sqrt{L})\|z\|_L,\hspB \|\lambda\|_{l^2}\leq \sigma_\max (A\sqrt{L})\|z\|_L, \label{bound DVM mom}
\end{gather}
where $\sigma_{\min/\max}(A\sqrt{L})$ represent the minimum/maximum singular value of the matrix $A\sqrt{L}$.  We will use the above two bounds to convert stability results for the DVM to stability results for the moment approximation.

\subsubsection{Collision step}
We start with the collision step given in \eqref{col step}. Applying triangle's inequality to the collision step we find 
\begin{equation}
\begin{aligned}
\mcal E_{k^*} \leq \frac{2}{(1+\Delta t/\tau)^2}\mcal E_{k} + 2\left(\frac{\Delta t/\tau}{1+\Delta t/\tau}\right)^2\sum_{i=1}^{N_x}\underbrace{\|\lan P_Mf_{\mcal M,i}^{k}\ran_N\|_{l^2}^2}_{\leq \sigma_{\max}(A\sqrt{L})^2\|W(f^k_{\mcal M,i})\|^2_{L}}.
 \label{stable col1}
\end{aligned}
\end{equation}
The bound on the right hand side follows from the inequalities in \eqref{bound DVM mom}. 
From page-92 of \cite{ConvgDVM} we know that the solution to the entropy-minimization problem \eqref{entropy min} satisfies
\begin{gather}
\|W(f^k_{\mcal M,i})\|^2_{L}\leq N^{3} \exp(2Nt_k). \label{stable maxwell}
\end{gather}
The above relation and the bound on $\mcal E_{k^*}$ given in \eqref{stable col1} provides 
\begin{gather}
\mcal E_{k^*} \leq \frac{2}{(1+\Delta t/\tau)^2}\mcal E_{k} + 2\left(\frac{\Delta t/\tau}{1+\Delta t/\tau}\right)^2N_x\sigma_{\max}(A\sqrt{L})^2N^{3} \exp(2Nt_k). 
\end{gather}
\subsubsection{Transport step}
With the following three steps, we establish the stability of the transport step given in \eqref{transport}. (i) We recover a DVM underlying the transport step in \eqref{transport}. (ii) Using stability properties of an upwind scheme, we establish the stability of the DVM. (iii) Finally, relating the discrete velocity solution to the moment solution, we establish the stability of the moment scheme. The details of these three steps is as follows.

We consider the reformulated transport step given in \eqref{conv step reform}. Let $\mcal N(AL)$ represent the null-space of the matrix $AL$, where $A$ and $L$ are as given in \eqref{mat-vec form} and \eqref{def A}, respectively. Then, the transport step provides 
\begin{equation}
\begin{aligned}
W(f_i^{k+1}) = &(1-\Lambda|\Xi|)W( f_i^{k^*})\\
& + \frac{\Lambda}{2}(|\Xi|-\Xi)W( f_{i+1}^{k^*})+\frac{\Lambda}{2}(|\Xi| + \Xi)W(f_{i-1}^{k^*}) + v,\label{conv step reform2}
\end{aligned}
\end{equation}
where $v$ belongs to $\mcal N(AL)$. Since the moments at time step $t_{k+1}$---given as $ALW(f_i^{k+1})$---are invariant under the choice of $v$, we choose $v = 0$. This makes the above evolution equation a space-time discretization of a system of decoupled linear advection equations given as $\pd_t W(f) + \Xi \pd_x W(f) = 0$. The discretization uses an explicit Euler and an upwind FV scheme to discretize the space and the time domain, respectively. From Example-7.2 of \cite{tadmorReview} we know that such a discretization is L2-stable under the CFL condition
\begin{gather}
0 < \Lambda \leq \opn{min}\{|\xiMax^{-1}|,|\xiMin^{-1}|\}/2. \label{stable lambda}
\end{gather}
This provides
\begin{gather}
\sum_{i=1}^{N_x}\|W(f_i^{k+1})\|^2_{L} \leq \sum_{i=1}^{N_x}\|W(f_i^{k^*})\|^2_{L} + |f_{in}(\cdot,t_k,\cdot)|^2_{\pd\Omega,N}.
\end{gather}
Above, $|\cdot|_{\pd\Omega,N}$ is as defined in \eqref{def bc norm}.
Using the bounds in \eqref{bound DVM mom} , we express the above bound in terms of moments to find
\begin{gather}
\mcal E_{k+1} \leq \kappa(A\sqrt{L})^2 \mcal E_{k^*} + \sigma_\max(A\sqrt{L})^2 |f_{in}(\cdot,t_k,\cdot) |^2_{\pd\Omega,N}.\label{bound transport}
\end{gather}
Above, $\kappa(A\sqrt{L})$ represents the condition number of the matrix $A\sqrt{L}$. We collect our stability estimate in the result below. 
\begin{theorem}
Consider the evolution scheme given in \Cref{sec: evo scheme} and let $\mcal E_k$ be the L2 energy defined in \eqref{def EM}. Assume that the boundary data satisfies \eqref{assume data} and that the ratio $\Lambda = \Delta t/\Delta x$ satisfies
$\Lambda\in (0,\opn{min}\{|\xiMax^{-1}|,|\xiMin^{-1}|\}/2].$ Then, $\mcal E_{k+1}$ is bounded as 
 \begin{gather}
 \mcal E_{k+1}\leq \mcal B_k + \mcal B_{\mcal M} + \mcal B_{in} \label{bound EM}
 \end{gather}
 where 
 \begin{equation}
 \begin{aligned}
 \mcal B_k := &\kappa(A\sqrt{L})^2 \frac{2}{(1+\Delta t/\tau)^2}\mcal E_{k}, \\
  \mcal B_{\mcal M}:= &2\kappa(A\sqrt{L})^2\sigma_{\max}(A\sqrt{L})^2\left(\frac{\Delta t/\tau}{1+\Delta t/\tau}\right)^2N_xN^{3} \exp(2Nt_k),\\
 \mcal B_{in} := &\sigma_\max(A\sqrt{L})^2 |f_{in}(\cdot,t_k,\cdot) |_{\pd\Omega,N}^2.
 \end{aligned}
 \end{equation}
\end{theorem}
We make the following remarks related to the above theorem.
\begin{enumerate}
\item The terms $\mcal B_k$, $\mcal B_{\mcal M}$ and $\mcal B_{in}$ appearing in \eqref{bound EM} represent the contribution from the previous time step, the discrete Maxwell-Boltzmann distribution function and the boundary data, respectively, into bound for the L2-energy at time $t_{k+1}$. Note that out of all these three terms, only $\mcal B_k$ depends upon the solution of the previous time-step.
\item For the limit $\tau \to 0$, at least formally, the BE results in the Euler equations \cite{chapman1990}. Under this limit, the bound in \eqref{bound EM} is robust, which is a result of performing the collision step implicitly.
\item The DVM corresponding to the transport step given in \eqref{conv step reform2} is a space-time discretization of a linear hyperbolic PDE. As a result, the L2-bound for the transport step (given in \eqref{bound transport}) is linear in time. In contrast, since the collision operator is non-linear, the collision step is non-linear. This introduces an exponential-in-time growth in the term $\mcal B_{\mcal M}$. 
\item For a fixed truncated velocity domain $\Omega_\xi$, consider the limit $N,M\to\infty$ with $N > M$. Under this limit, the bound on $\mcal E_{k+1}$ is not robust because---at least heuristically---both $\kappa(A\sqrt{L})$ and $\sigma_\max(A\sqrt{L})$ are almost independent of $N$ and grow polynomially with $M$. To derive bounds that are independent of $\kappa(A\sqrt{L})$ and $\sigma_\max(A\sqrt{L})$, one should directly consider the moment approximation without accessing the underlying DVM. As yet, it is unclear how to proceed with such a technique. 
\item Nowhere in the proof of the above theorem we used the fact that we minimize the L2-norm in the moment-closure problem given in \eqref{l2 quad}. Therefore, the bound on $\mcal E_{k+1}$ holds for any other objective functional and specifically for the minimum-entropy closure considered in \cite{GrothMaxEntropy,Roman35}. To the best of our knowledge, none of the other works develop such a bound for a minimization-based closure.  
\end{enumerate}

\subsection{Computational costs} We study the cost of evolution scheme outlined in \Cref{sec: evo scheme}. We consider the cost of a single time-step performed in a single spatial cell.
\begin{enumerate}
\item \textit{Entropy-minimization step:} We use Newton iteration to solve the entropy-minimization problem where we compute and invert a Hessian $H(x,t)$ given as
\begin{gather}
\left(H(x,t)\right)_{kl} := \sum_{i}\left( \Pcons(\xi_i)\right)_k\left( \Pcons(\xi_i)\right)_l \exp(\Pcons\cdot \alpha(x,t))\omega_i. \label{hessian entropy}
\end{gather}
 Computing the Hessian is an $\mcal O(N)$ operation. As a stopping criterion to the Newton solver, we consider a user-defined tolerance of \texttt{TOL} in the moment constraints. Suppose we need $m_{\texttt{TOL}}$ Newton iterations to reach this tolerance then, the total cost of entropy-minimization is given as 
$$
C_{\opn{entropy}} = \mcal O(N m_{\texttt{TOL}}).
$$
In all our numerical examples, we choose $\texttt{TOL} = 1e-8$.
\item \textit{Collision step:} Computing the $M$-moments of the discrete Maxwell-Boltzmann pdf is an $\mcal O(N M)$ operation and updating the moments in the collision step is an $\mcal O(M)$ operation. Thus, the total cost of the collision step is given as 
$$
C_{\opn{col}} = \mcal O(MN).
$$
\item \textit{Optimization step:} We use the \texttt{quadprog} routine from matlab to solve the optimization problem in \eqref{l2 quad} and we use the default interior-point-convex solver with all the parameters set to their default values. Usually, it is difficult to estimate the complexity of this algorithm but a crude estimate gives \cite{ComplexIntPoint}
\begin{gather}
C_{\opn{opt}} =\mcal O(N^3).
\end{gather}
\item \textit{Transport step:} Flux computation is an $\mcal O(MN)$ operation and the time update of the moments is an $\mcal O(M)$ operation. Thus the cost of the transport step is 
$$
C_{\opn{tran}} = \mcal O(MN).
$$
\end{enumerate}
Summing up the above costs, the total cost of our evolution scheme is given as
$$
C_{\opn{total}} = \mcal O(N m_{\texttt{TOL}}) + \mcal O(MN) + \mcal O(N^3).
$$

\begin{remark}[Efficiency of the optimization step]
For $N > M$ (the values of $N$ that interest us, see \Cref{remark: value N}) and a sufficiently small $m_{\texttt{TOL}}$, solving the quadratic optimization problem is the most expensive part of the algorithm. There are two possible way to reduce this cost (i) choose $M,N$ adaptively and vary them over the space-time domain \cite{MalikAdaptivity,Torrilhon2017}, or (ii) train an auto-encoder/gaussian-regression to solve the optimization problem \cite{GPEntropy,Han2019}. We plan to consider both these directions in the future. 
\end{remark}

\section{Extension to multi-dimensions}\label{sec: extd multi}
Maintaining consistency with our numerical experiments, we propose an extension of our method to two-dimensional planar flows. An extension to three-dimensional problems is similar and is not discussed for brevity. For 2D problems, we reduce the storage requirements by solving for the reduced pdfs $h_1$ and $h_2$ given as \cite{PlanarFlows}
\begin{equation}
\begin{aligned}
h_1(x,t,\xi_1,\xi_2) := &\int_{\mbb R}f(x,t,\xi_1,\xi_2,\xi_3)d\xi_3,\\
 h_2(x,t,\xi_1,\xi_2) := &\int_{\mbb R}\xi_3^2f(x,t,\xi_1,\xi_2,\xi_3)d\xi_3.
\end{aligned}
\end{equation}
In the coming discussion, $\xi$ will represent a velocity vector in $\mbb R^2$ and with $\lan g\ran$ we will represent the integral of a function $\xi\mapsto g(\xi)$ over $\mbb R^2$. To derive the governing equation for $h_1$ and $h_2$, we multiply the BTE given in \eqref{BTE} by $1$ and $\xi_3^2$ and integrate over $\mbb R$ with respect to $\xi_3$ to find
\begin{gather}
\pd_t h_i + \xi_1 \pd_{x_1}h_i + \xi_2 \pd_{x_2}h_i = \frac{1}{\tau(x,t)}\left(h_{i,\mcal M}-h_i\right).\label{BTE 2D}
\end{gather}
Above, $h_{i,\mcal M}$ represents the reduced Maxwell-Boltzmann pdf and is given as 
\begin{gather}
h_{1,\mcal M} = \frac{\rho}{2\pi\theta}\exp\left(-\frac{|\xi-v|^2}{2\theta}\right),\hspB h_{2,\mcal M} = \frac{\rho}{2\pi} \exp\left(-\frac{|\xi-v|^2}{2\theta}\right),
\end{gather}
where, $|\cdot|$ is the Eucledian norm of a vector. Note that the mass $\rho$, the momentum $\rho v$ and the temperature $\theta$ can be recovered from $h_1$ and $h_2$ via 
\begin{equation}
\begin{gathered}
\rho = \lan h_1 \ran,\hspB \rho v = \lan \xi h_1\ran,\hspB \rho\theta = \frac{1}{3}\left(\lan |\xi|^2h_1\ran - \rho|v|^2 + \lan h_2\ran \right).
\end{gathered}
\end{equation}

\subsection{Moment equations}
The moment approximation we discuss below is the same for both $h_1$ and $h_2$. Therefore, for the simplicity of notation, we present our approximation for some representative $h$. Similar to the 1D case, we truncate the velocity domain to $\mbb R^2\supset\Omega_\xi = [\xi_{1,\min},\xi_{1,\max}]\times [\xi_{2,\min},\xi_{2,\max}]$. To compute $\xi_{i,\max/\min}$, we adopt the same methodology as that outlined in \Cref{sec: compute xiMax}. We consider tensorized $N\times N$ Gauss-Legendre quadrature points inside $\Omega_\xi$. Using these quadrature points, we approximate $\lan\cdot\ran$ by $\lan\cdot\ran_{N,N}$.

To derive a governing equation for the moments of $h$, we first define a polynomial in $\xi$. With $\beta_M:=(\beta^M_1,\beta^M_2)\in\mbb R^2$ we represent a multi-index with each entry being a natural number and the $l^1$ norm of $\beta_M$ being equal to $M$. Using $\beta_M$, we define a $M$-th order polynomial in $\xi$ via $p_{\beta_M} = \xi_1^{\beta^M_1}\xi_2^{\beta^M_1}$. Note that for a given $M$, $\beta_M$ is non-unique---for $M=1$, $\beta_M$ could either be $(0,1)$ or $(1,0)$. In a vector $P_M(\xi)$, we collect all the polynomials $p_{\beta_M}$ upto order $M-1$. For completeness, we present the entries in $P_M(\xi)$ for $M=3$ and $M=5$.
\begin{equation}
\begin{aligned}
\text{M=3: }P_M(\xi)= &\left(1,\xi_1,\xi_2,\xi_1^2,\xi_1\xi_2,\xi_2^2\right)^T;\\
\text{M=5: }P_M(\xi)= &\left(1,\xi_1,\xi_2,\xi_1^2,\xi_1\xi_2,\xi_2^2,\xi_1^3,\right.\\
&\left. \xi_1^2\xi_2,\xi_1\xi_2^2,\xi_2^3,\xi_1^4,\xi_1^3\xi_2^1,\xi_1^2\xi_2^2,\xi_1\xi_2^3,\xi_2^4\right)^T.
\end{aligned}
\end{equation}
Note that for $M=3$ and $M=3$, $P_M(\xi)$ has $6$ and $15$ entries, respectively.

For some $M\in \mbb N$, we approximate $h$ by $h_M$ where we compute $h_M$ using the L2-minimization problem given in \eqref{l2 quad}. To evolve the moments of $h_M$, we use a multi-dimensional version of the moment system given in \eqref{mom 1D}, which reads 
\begin{equation}
\begin{aligned}
\pd_t\lan P_M\hMN \ran_{N,N}+ &\pd_{x_1} \lan P_M \xi_1\hMN \ran_{N,N} + \pd_{x_2}\lan P_M \xi_2\hMN \ran_{N,N}\\
&= \frac{1}{\tau}(\lan P_M h_{\mcal M,N}\ran_{N,N}-\lan P_M h_{M}\ran_{N,N})\hsp\text{on}\hsp \Omega\times D,\\
\quad\hsp \lan P_M\hMN(t=0)\ran_{N,N}=&\lan P_M h_{0}\ran_{N,N}\hsp\text{on}\hsp\Omega. \label{mom 2D}
\end{aligned}
\end{equation}
Above, $h_{\mcal M,N}$ is a discretization of the Maxwell-Boltzmann pdf that results from solving a multi-dimensional version of the optimization problem given in \eqref{entropy min}---see \cite{LucDVM} for an explicit form of this optimization problem. The treatment of boundary conditions is the same as that for the 1D case and is not discussed for brevity. 
\subsection{Space-time discretization}
For simplicity, we consider a square spatial domain $\Omega = [x_{1,\min},x_{1,\max}]\times [x_{2,\min},x_{2,\max}]$. We discretize $\Omega$ with $N_x$ number of uniform elements in each spatial dimension and with $\Delta x$ we represent the grid spacing. With some additional technical details, it is straightforward to extend our framework to curved domain discretized with unstructured meshes. For simplicity, we consider a fixed time-step of size $\Delta t$.

We index a spatial cell with $(i,j)$ where $i,j\in \{1,\dots,N_x\}$. With $\lan P_M h_{i,j}^k \ran_{N,N}$ we represent a FV approximation to $\lan P_M \hMN(x,t_k,\cdot) \ran_{N,N}$ in the cell $\mcal I_{i,j}$. Given $\lan P_M h_{i,j}^k \ran_{N,N}$, we want to compute the FV approximation at the next time instance. To this end, we follow the same four steps as those outlined for the 1D-case in \Cref{sec: evo scheme}. The entropy-minimization step, the collision step and the optimization step are very similar to the 1D case and, for brevity, we do not repeat them here. The transport step is slightly different and is given as 
\begin{equation}
\begin{aligned}
\frac{\lan P_M f_{i,j}^{k+1}\ran_{N,N}-\lan P_M f_{i,j}^{k^*}\ran_{N,N}}{\Delta t}= &-\frac{1}{\Delta x}\left(\mcal F_1(W(f^{k,*}_{i+1,j}),W(f^{k,*}_{i,j}))\right.\\
&\left.-\mcal F(W(f^{k,*}_{i,j}),W(f^{k,*}_{i-1,j}))\right)\\
&-\frac{1}{\Delta x}\left(\mcal F_2(W(f^{k,*}_{i,j+1}),W(f^{k,*}_{i,j}))\right.\\
&\left.-\mcal F_2(W(f^{k,*}_{i,j}),W(f^{k,*}_{i,j-1}))\right).
\end{aligned}
\end{equation}
Above, $\{W(f^{k,*}_{i,j})\}_{i,j}$ results from the optimization step and $\mcal F_1(W_1,W_2)$ and $\mcal F_2(W_1,W_2)$ are the numerical fluxes given as 
\begin{gather}
\mcal F_i(W_1,W_2) :=\frac{1}{2}\left(AL(\Xi_i-|\Xi_i|)W_1 + AL(\Xi_i+|\Xi_i|)W_2 \right).
\end{gather}
Above, $A$ and $L$ are multi-dimensional versions of the matrices given in \eqref{mat-vec form} and $\Xi_i$ is a diagonal matrix with all the $i$-th components of the quadrature point's locations at its diagonal. 

Assuming that the initial and the boundary data satisfies \eqref{assume data}, one can show that the space-time discretization results in a feasible optimization if the ratio $\Lambda = \Delta t/\Delta x$ satisfies 
\begin{gather}
0 < \Lambda \leq \frac{1}{2}\min_i\left\{\min\left\{|\xi_{i,\max}^{-1}|,|\xi_{i,\min}^{-1}|\right\}\right\}. 
\end{gather}
Similarly, one can show that the space-time discretization is L2-stable if $\Lambda$ satisfies
\begin{gather}
0 < \Lambda \leq \frac{1}{4}\min_i\left\{\min\left\{|\xi_{i,\max}^{-1}|,|\xi_{i,\min}^{-1}|\right\}\right\}. 
\end{gather}
A proof of the above two results uses the exact same technique as that for the 1D case and is not repeated for brevity.

\section{Numerical Results}\label{sec: num exp}
For simplicity, we non-dimensionalize the BE and all the macroscopic quantities with appropriate powers of some reference density $\rho_0$, temperature $\theta_0$ and length scale $l$. This introduces the Knudsen number $Kn$ that scales the collision operator $Q(f)$ and reads $Kn := \tau_0/\left(\sqrt{\theta_0}l\right)$---we refer to \cite{Struchtrupbook} for the details of non-dimensionalization. In the definition of the collision frequency $\tau(x,t)^{-1}$ given in \eqref{BGK op}, we choose $C=1$ and $\omega = 1$. Our choice of $C$ and $\omega$ does not necessarily corresponds to a physical system and is made for demonstration purposes. 

We consider the following test cases. 
\begin{enumerate}
\item \textbf{Test case-1} We consider the pdf
\begin{gather}
f(\xi) = \frac{1}{\sqrt{2\pi\theta_0}}\exp\left(-\frac{(\xi-u_0)^2}{2\theta_0}\right) +  \frac{1}{\sqrt{2\pi\theta_1}}\exp\left(-\frac{(\xi-u_1)^2}{2\theta_1}\right).\label{bi Gauss}
\end{gather}
Given the first $M$ moments of $f$ and using the pos-L2-MOM, we approximate the $M+1$-st moment of $f$. We study the error of this approximation with respect to the number of moments $M$. We choose $\theta_0 = 3$, $u_0 = -4$, $\theta_1 = 4$ and $u_1=5$, which ensures that $f$ is far away from a Maxwell-Boltzmann distribution function in the Kullback–Leibler divergence sense. 
\item \textbf{Test case-2} For a one-dimensional space-velocity domain, we consider the Sod's shock tube problem from \cite{R13Simulation}. We set $\Omega = [-2,2]$ and $D = [0,0.3]$. Recall that $D$ is the time domain. As the initial data, we consider a gas at rest and at equilibrium. We initialize the temperature $\theta$ with a constant value of one and we initialize density as 
\begin{gather}
\rho(x,t=0)=\begin{cases}
7,\hspB x\leq 0\\
1, \hspB x>0
\end{cases}.
\end{gather}
As the boundary data $f_{in}$, we consider a Maxwell-Boltzmann pdf. At $x=\xMin$ and for all $t\in D$, we set density to $7$, velocity to $0$ and temperature to $1$. The velocity and the temperature at the right boundary remains the same but the density changes to one. We consider two different values of the Knudsen number---$Kn=0.1$ and $Kn=0.01$.
\item \textbf{Test case-3} For a one-dimensional space-velocity domain, we consider the two-beam interaction experiment from \cite{Roman35}. The space-time domain $\Omega\times D$ remains the same as the previous test case. As the initial data, we consider a gas at equilibrium with a constant density and temperature of one. As the initial velocity, we consider 
\begin{gather}
v(x,t=0)=\begin{cases}
1,\hspB &x\leq 0\\
-1, \hspB &x>0
\end{cases}.
\end{gather}
As the boundary data $f_{in}$, we consider a Maxwell-Boltzmann pdf. At $x=\xMin$ and for all $t\in D$, we set density to $1$, velocity to $1$ and temperature to $1$. The density and the temperature at the right boundary remains the same but the velocity changes to $-1$. We consider two different values of the Knudsen number---$Kn=0.1$ and $Kn=0.01$.
\item \textbf{Test case-4} We consider a two-dimensional spatial domain and a planar flow regime. We choose $\Omega = [0,1]\times [0,1]$ and $D=[0,0.2]$. We consider a micro-bubble dispersion problem where we start with a fluid at equilibrium and at rest. We consider a constant temperature of one and consider a density given as
 \begin{gather}
\rho(x,t=0)=\rho_0 + \exp(-|x-1|^2\times 10^{2}),\hspB\forall x\in\Omega.
\end{gather}
As the ground state density, we set $\rho_0=1$. As the boundary data $f_{in}$, we consider a Maxwell-Boltzmann pdf with a density $\rho_0$, velocity zero and temperature one. We consider a Knudsen number of $0.1$. 

We emphasis that for this test case, it is crucial that the moment-closure problem has a positive solution. Otherwise, the density can get negative resulting in a breakdown of the solution algorithm. We refer to \cite{SarnaHarshit2020} for a similar experiment involving the linearized BE and the Grad's Hermite expansion, which is not necessarily positive. There, the deviation in density gets negative for small values of $M$. However, since the BE is linearized, negative densities do not crash the solution algorithm. 
\end{enumerate}
\subsection{Test case-1} We truncate the velocity domain to $\Omega_\xi = [-20,20]$. This ensures that the support of $f$ (upto machine precision) is contained inside $\Omega_\xi$. We compute $f_M$ using the optimization problem given in \eqref{l2 quad}. 

\subsubsection{Error in the higher order moment} Recall that we used $f_M$ to close the moment system in \eqref{mom 1D} by approximating the $M$-th order moment of $f$. The relative error of this approximation is given as
\begin{gather}
\mcal E(M):=\left|\frac{\lan \xi^{M} (\fMN-f)\ran_N}{\lan \xi^{M} f\ran_N}\right|.\label{def EM}
\end{gather}
We study $\mcal E(M)$ for different values of $M$. We vary $M$ from $3$ to $22$ in steps of one, and we fix $N$ at a sufficiently large value of $40$.

As $M$ increases, $\mcal E(M)$ appears to converge to zero, although not monotonically---see \Cref{fig: test1 error}. Note that this non-monotonic convergence is typical also for a Grad's moment approximation \cite{Torrilhon2015,SarnaHarshit2020,Cai2017}. However, unlike the Grad's moment approximation where the error convergences monotonically for either the even or the odd values of $M$, the convergence behaviour of the pos-L2-MOM is rather random. For instance, the error (slightly) increases from $M=5$ to $M=7$. Similarly, the error (slightly) increases from $M=15$ to $M=17$. Note that for $M\geq 16$, the error appears to converge monotonically. 

\subsubsection{Error in approximating the pdf}
\begin{figure}[ht!]
\centering
\subfloat[Convergence behaviour ]{\label{fig: test1 error}\includegraphics[width=2.8in]{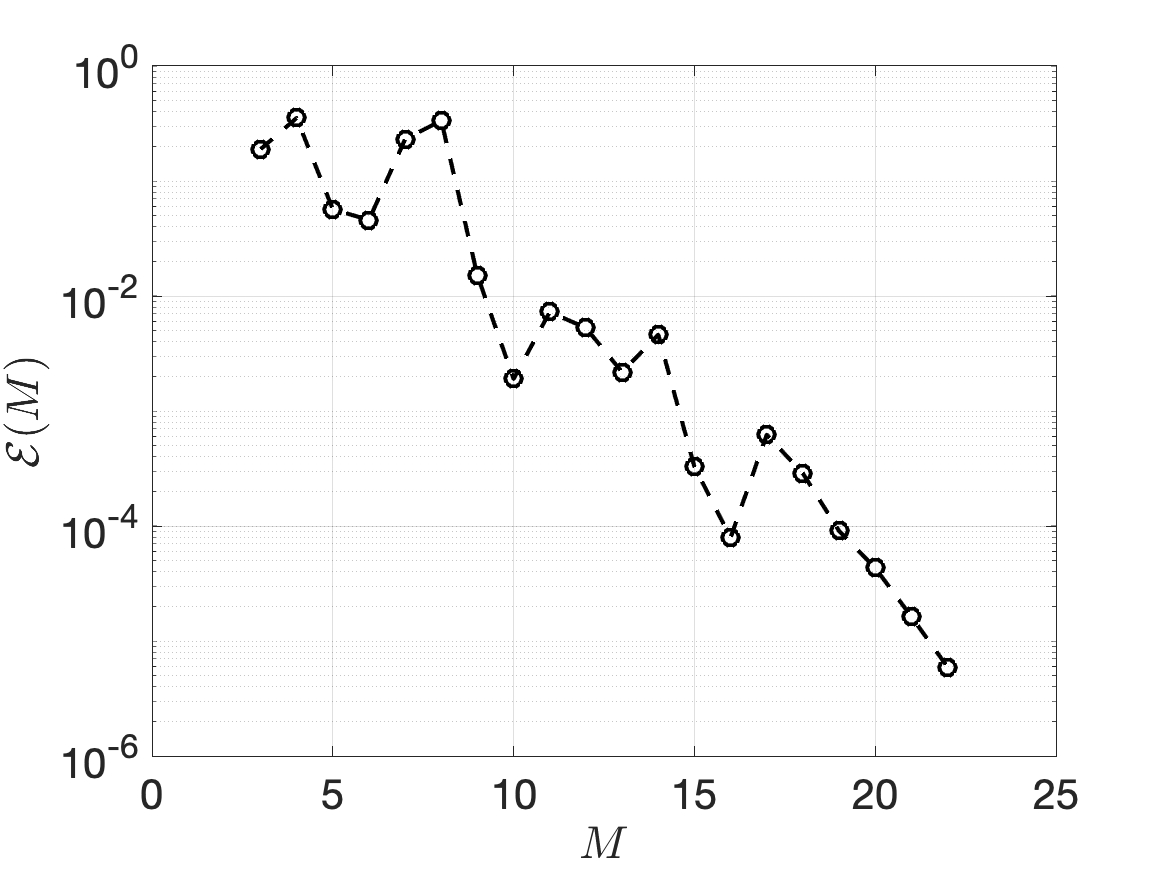}}
\subfloat[$f$ and $f_M$]{\label{fig: test1 f}\includegraphics[width=2.8in]{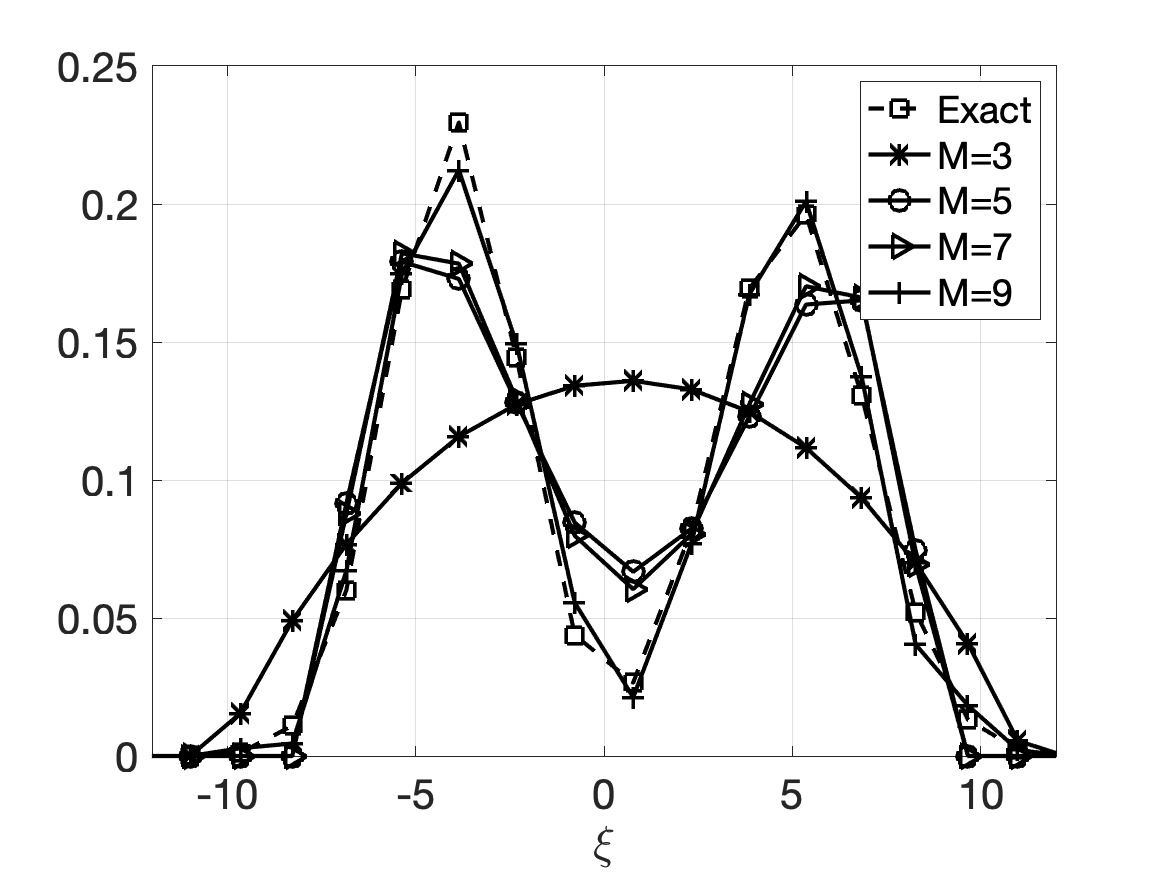}}
\hfill
\subfloat[$f$ and $f^{DG}_M$]{\label{fig: test1 fDG}\includegraphics[width=2.8in]{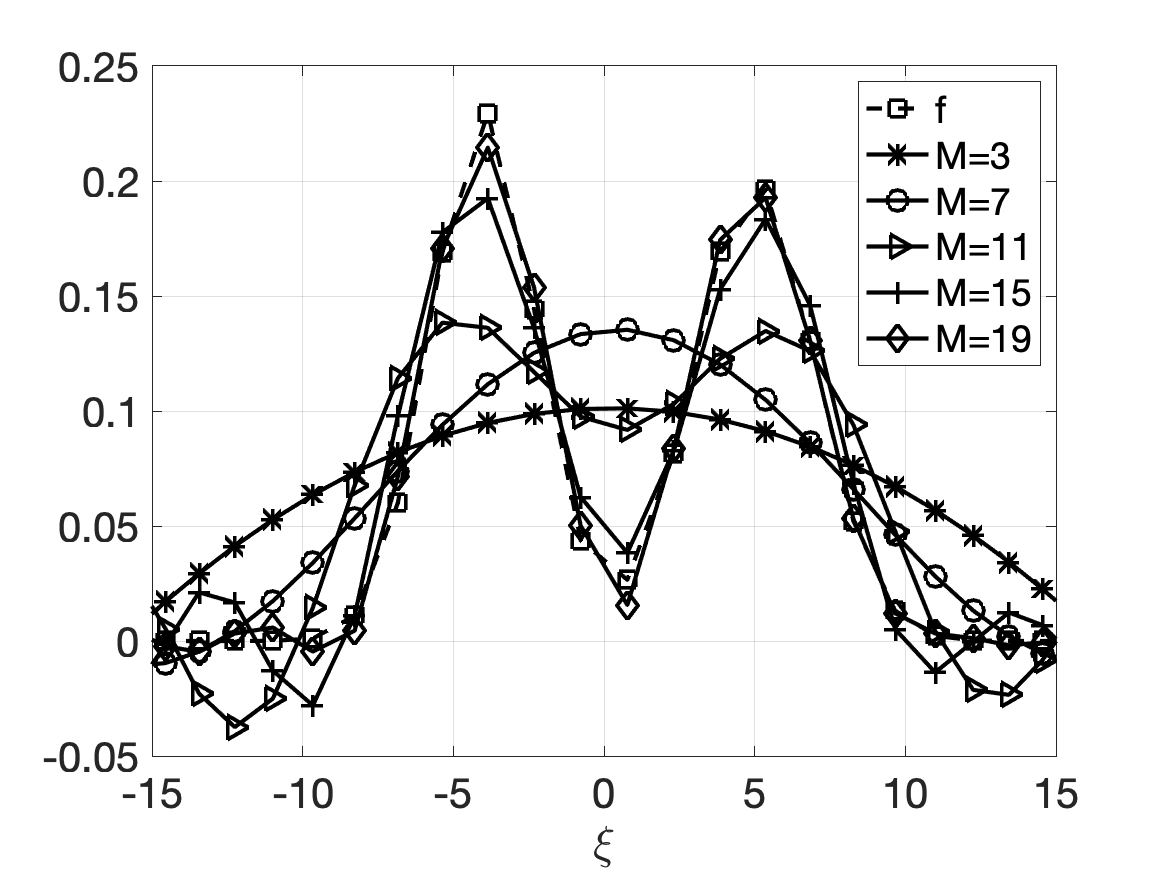}}
\subfloat[Comparison of L2-error]{\label{fig: test1 comp_DG}\includegraphics[width=2.8in]{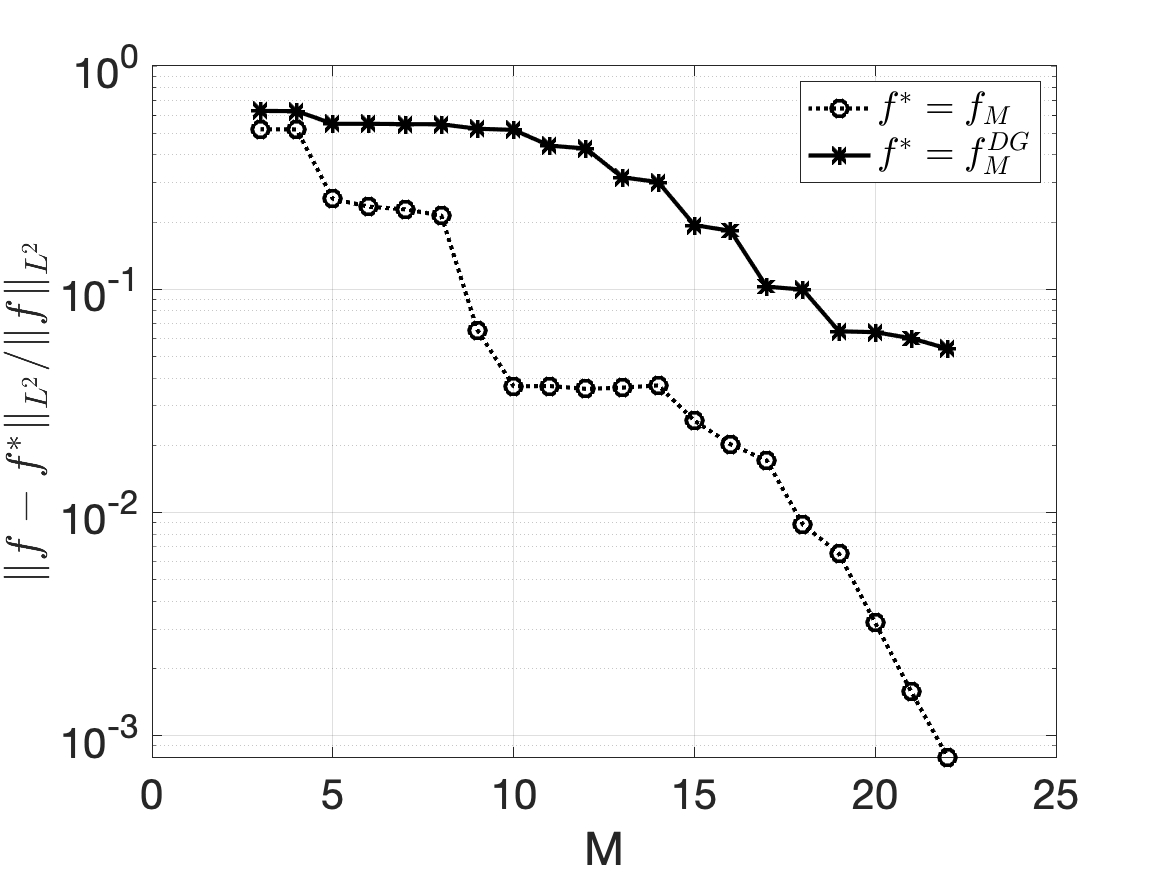}}
\hfill
\caption{Results for test case-1. (a) and (b) y-axis is on a log-scale.} 
\end{figure}

For different values of $M$, \Cref{fig: test1 f} compares $f$ to $f_M$. To extend the discrete values of $f_M$ to $\Omega_\xi$, we perform a piecewise linear interpolation between the quadrature points. For $M=3$, pos-L2-MOM is unable to capture the general shape of the function. Nevertheless, increasing the value of $M$ improves the results. Already for $M=5$, we observe that $f_M$ has two distinct peaks and starts to capture the shape of the function. Increasing $M$ from $5$ to $7$ does not show much of an improvement. However, increasing $M$ from $7$ to $9$ improves the results significantly. The result for $M=9$ almost overlaps the exact solution with little deviations. Let us mention that for all values of $M$, $f_M$ remains positive. 

For a comparison, we compute a DG approximation of $f$. We represent the DG approximation by $f_M^{DG}$ and compute it by projecting $f$ (under the $L^2(\Omega_\xi)$ inner-product) onto the first $M$ Legendre polynomials in $\xi$. For the different values of $M$, \Cref{fig: test1 fDG} compares $f$ to $f_{M}^{DG}$. Since a DG approximation does not penalize negativity (see \Cref{remark: DG}), for all values of $M$, $f_M^{DG}$ is negative for some part of the velocity domain. Furthermore, only for $M\geq 11$, the DG approximation starts to capture the general shape of the function. Compare this to $f_M$, which, already for $M=5$, accurately represents the shape of the function. 

The superior accuracy of $f_M$---as compared to $f_M^{DG}$---in approximating $f$ is clearly visible in \Cref{fig: test1 comp_DG}, which compares the relative L2-error in approximating $f$. The difference between the error values becomes larger as the value of $M$ increases. For the largest value of $M$ equals $22$, the relative L2-error resulting from the approximation $f_M$ is $8\times 10^{-4}$, which is $\approx 10^{-2}$ times smaller than that resulting from the approximation $f_M^{DG}$. 

\subsection{Test case-2}
\subsubsection{Reference solution}\label{algo ref sol}
We compute the reference solution using a DVM proposed in \cite{LucDVM}. We consider an explicit Euler time-stepping scheme and a first-order FV spatial discretization. We truncate the velocity domain to $[-7,7]$, and place $N  = 350$ velocity grid points inside the truncated velocity domain. As the velocity grid points, we consider Gauss-Legendre quadrature nodes. We discrete the space domain with $N_x=10^3$ uniform cells and consider a constant time-step of $\Delta t = 0.5 \times \Delta x/7$. To arrive at these discretization parameters, we performed a convergence study that consisted of the following steps. (i) With the velocity and the temperature field taken from the initial data, estimate $\xiMaxMin$ using the relation in \eqref{def ximax}. For the present test case, this provides $\xiMax = 3.5$ and $\xiMin = -3.5$. (ii) Fix $N_x$ at $10^3$ and $\Delta t$ to $0.5\times \Delta x/\xiMax $.  (iii) Choose $N = 50$ and increase it to $350$ in steps of $50$. (iv) Terminate the refinement as soon as the relative change in mass, momentum and energy between two subsequent refinement cycles drops below a tolerance of $10^{-5}$. (v) If the tolerance is not reached, increase $\xiMax$ by $0.5$, decrease $\xiMin$ by $0.5$ and repeat the process from step-(ii). Note that if the refinement cycle does not terminate then one should increase the value of $N_x$ and repeat the entire process. For the all the earlier mentioned test cases, the value of $N_x=10^3$ was sufficiently large to terminate the refinement cycle. 
\subsubsection{Convergence study}
We are interested in the relative L2 error in the different macroscopic quantities that we define as 
\begin{gather}
\mcal E_{cons}(M,N_x):=\frac{\|\lan\Pcons (\fMNNx(\cdot,t=T,\cdot)-f_{DVM}(\cdot,t=T,\cdot))\ran_N\|_{L^2(\Omega;\mbb R^3)}}{\|\lan\Pcons f_{DVM}(\cdot,t=T,\cdot)\ran_N\|_{L^2(\Omega;\mbb R^3)}}. 
\end{gather}
Above, $\Pcons$ and $\fMNNx$ are as defined in \eqref{def Pcons}  and \eqref{FV fmn}, respectively. We keep the value of $N$ fixed at $30$. 

\begin{figure}[ht!]
\centering
\subfloat[]{\label{fig: test2 err Kn0p1}\includegraphics[width=2.8in]{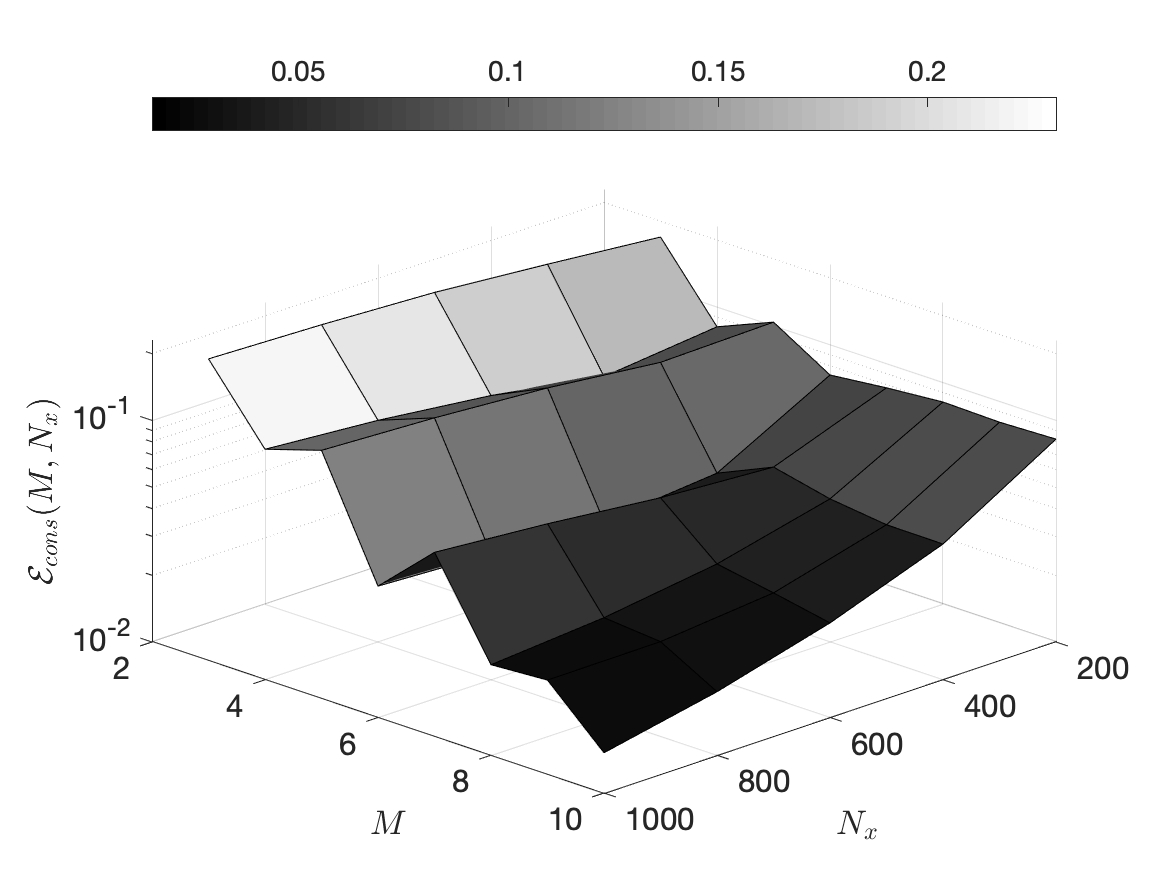}}
\subfloat[]{\label{fig: test2 err Kn0p01}\includegraphics[width=2.8in]{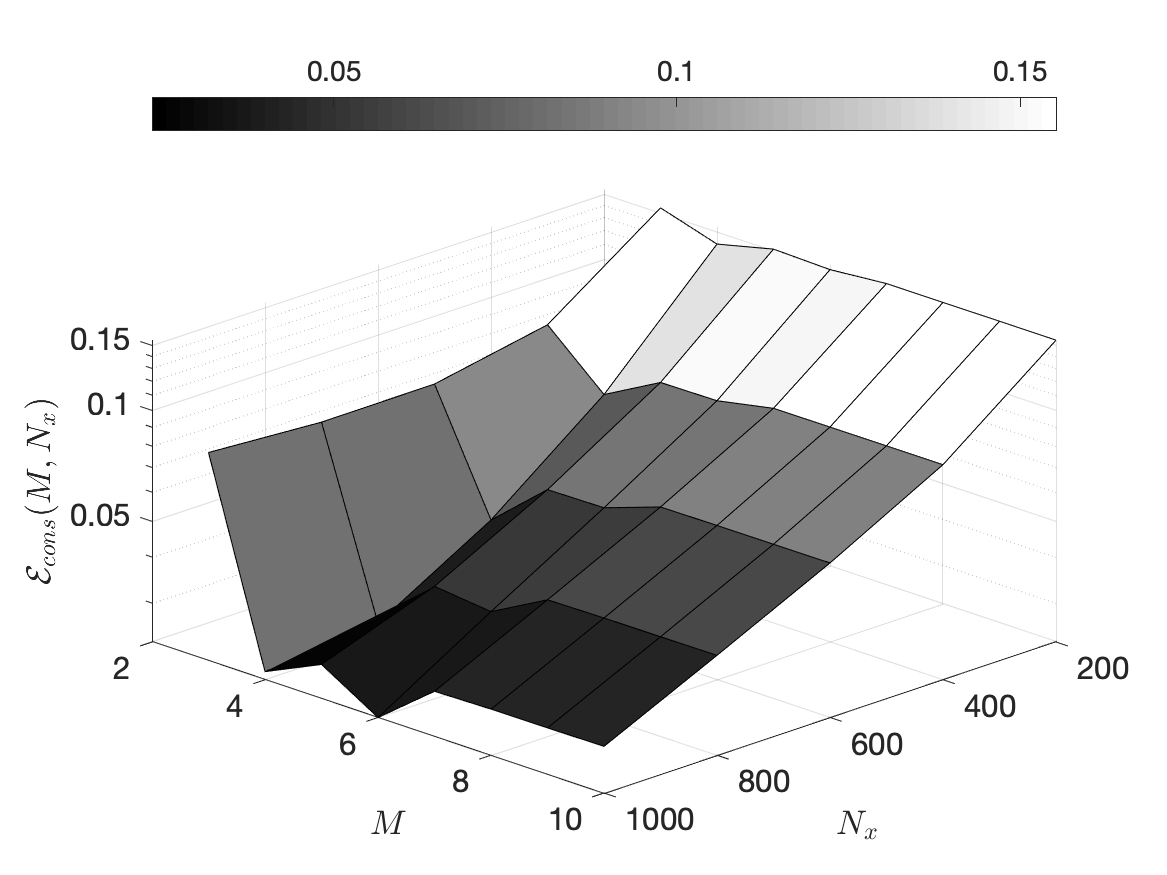}}
\caption{Results for test case-2. Convergence of the relative error with $N_x$ and $M$. Computations performed using the pos-L2-MOM. (a) $\opn{Kn}=0.1$ and (b) $\opn{Kn} = 0.01$. The z-axis on both the plots is on a log-scale.} 
\end{figure}

We first consider $\opn{Kn} =0.1$. We increase $M$ from $3$ to $10$ in steps of $1$ and $N_x$ from $200$ to $10^3$ in steps of $200$. We choose $\Delta t = 0.5 \times \Delta x/7$. \Cref{fig: test2 err Kn0p1} shows the error $\mcal E_{cons}(M,N_x)$ for the different values of $M$ and $N_x$. Fixing $N_x$ at a small value---$200$ for instance---and increasing $M$ does not reduce the error. This is because for small values of $N_x$, the error is dominated by the error in our space-time discretization. Furthermore, for a small value of $M$, increasing $N_x$ beyond a certain limit does not decrease the error. On the other hand, choosing a large value of $N_x$---$10^3$ for instance---and increasing $M$, or increasing both $M$ and $N_x$ simultaneously, reduces the error. Note that similar to the previous test case, the error decay is not monotonic. Our results suggest that to balance the accuracy with the computational cost, an adaptive choice of $M$ and an adaptive spatial grid is desirable. We plan to develop such an adaptive framework in the future---see \cite{MalikAdaptivity} for an adaptive moment method. Let us also mention that at $N_x=10^3$ and $M=10$, we attain a minimum relative error of $2.4\times 10^{-2}$. We find this error value acceptable, given that $M=10$ is less than $10\%$ of the velocity grid points used in our reference DVM. 

We now consider $\opn{Kn} = 0.01$. We choose $M$ and $N_x$ as before. \Cref{fig: test2 err Kn0p1} shows the error $\mcal E_{cons}(M,N_x)$ for the different values of $M$ and $N_x$. As compared to $\opn{Kn} = 0.1$, the smaller values of $M$ perform much better, which is in accordance with similar studies conducted in the previous works \cite{Torrilhon2015}. For instance, consider the results for $M=4$ and $N_x = 10^3$. For $\opn{Kn} = 0.1$, we find $\mcal E_{cons}(4,10^3) = 1.3\times 10^{-1}$, whereas for $\opn{Kn}=0.01$ we find $\mcal E_{cons}(4,10^3) = 2.5\times 10^{-2}$, which is almost an order-of-magnitude better than the result for $\opn{Kn} = 0.1$. 

Although the lower values of $M$ perform better for $\opn{Kn}=0.01$ than for $\opn{Kn}=0.1$, the minimum error attained is almost the same for both the Knudsen number---for $\opn{Kn}=0.01$ the minimum error is $2.3\times 10^{-3}$, which is $0.95$ times that of the minimum error for $\opn{Kn} = 0.1$. This is because for $\opn{Kn}=0.01$, the error at $N_x=10^3$ is already dominated by the error in our spatial discretization and we see almost no error reduction upon increasing $M$ from $7$ to $10$. By increasing $N_x$ from $10^3$ to $1.5\times 10^3$, we could remove this error stagnation and for $M=10$, achieve an error of $1.2\times 10^{-3}$. 
\subsubsection{Sub-shocks}
\begin{figure}[ht!]
\centering
\includegraphics[width=2.8in]{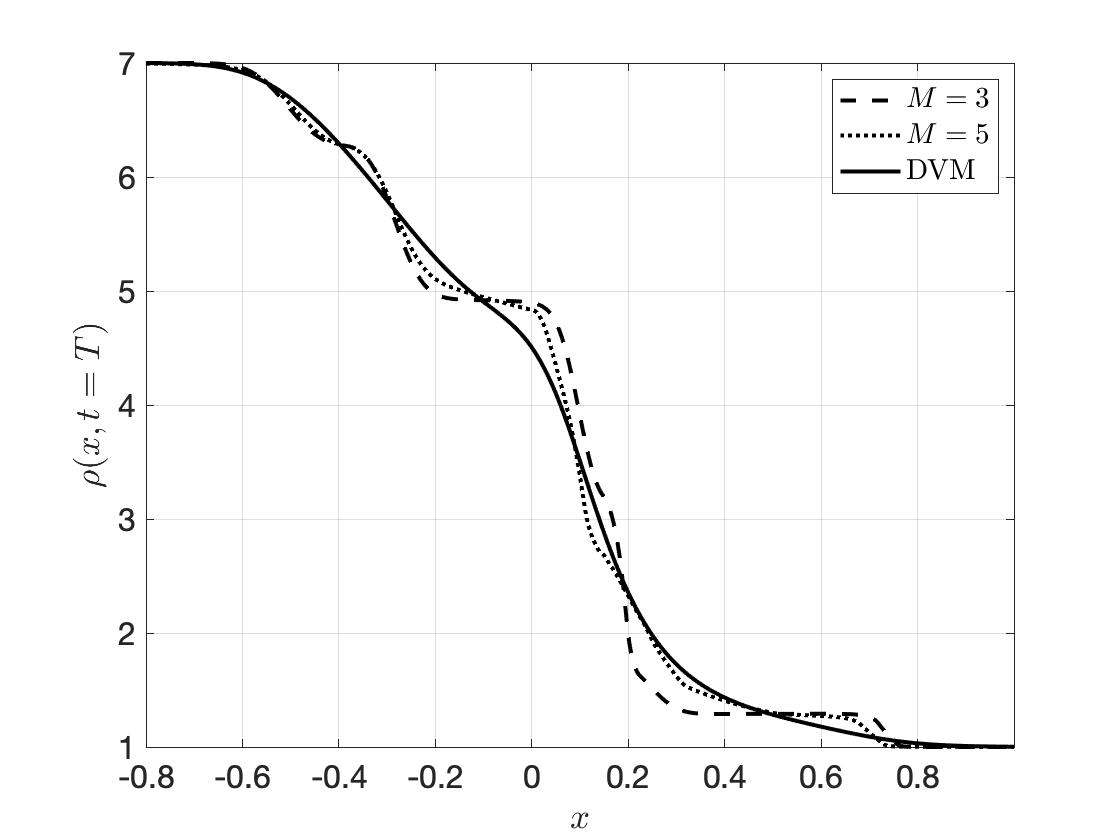}
\caption{Results for test case-2. Density profile for $Kn=0.1$ and at $t=T$. Computations performed with $N_x=10^3$ grid-cells.} 
\label{fig: test2 sub-shock}
\end{figure}

Shock speeds that are faster than the characteristic speeds in a moment system result in sub-shocks---we refer to \cite{Torrilhon2000} for an exhaustive study of sub-shocks for the Grad's MOM. Similar to the Grad's MOM, the pos-L2-MOM shows sub-shocks-type structures---see the density profile shown in \Cref{fig: test2 sub-shock}. These structures have a staircase-type shape, and increasing $M$ from $3$ to $5$ has a smoothing effect that reduces the staircase effect. To conclude that these structures are indeed sub-shocks, one needs to study the characteristic speeds of the moment system given in \eqref{mom 1D}. Note that these sub-shocks can be removed by introducing second-order spatial derivatives in the moment equations via regularization---see the discussion on the regularized-13 moment equations \cite{ShockStructureR13}.

\subsection{Test case-3}
As before, we construct a reference solution using the DVM. The convergence study discussed in \Cref{algo ref sol} lead to $N_x=10^3$, $\xiMax = 5$, $\xiMin = -5$ and $N = 350$. For the pos-L2-MOM, we fix $N = 30$ and $N_x=10^3$, and study the results for two different values of $M$, $M=5$ and $M=7$. We choose $\Delta t = 0.5 \Delta x/\xiMax$. The convergence behaviour is similar to the previous test case and not discussed for brevity.

\begin{figure}[ht!]
\centering
\subfloat[$M=5$, $\opn{Kn}=0.1$]{\label{fig: test3 M5_Kn0p1}\includegraphics[width=2.8in]{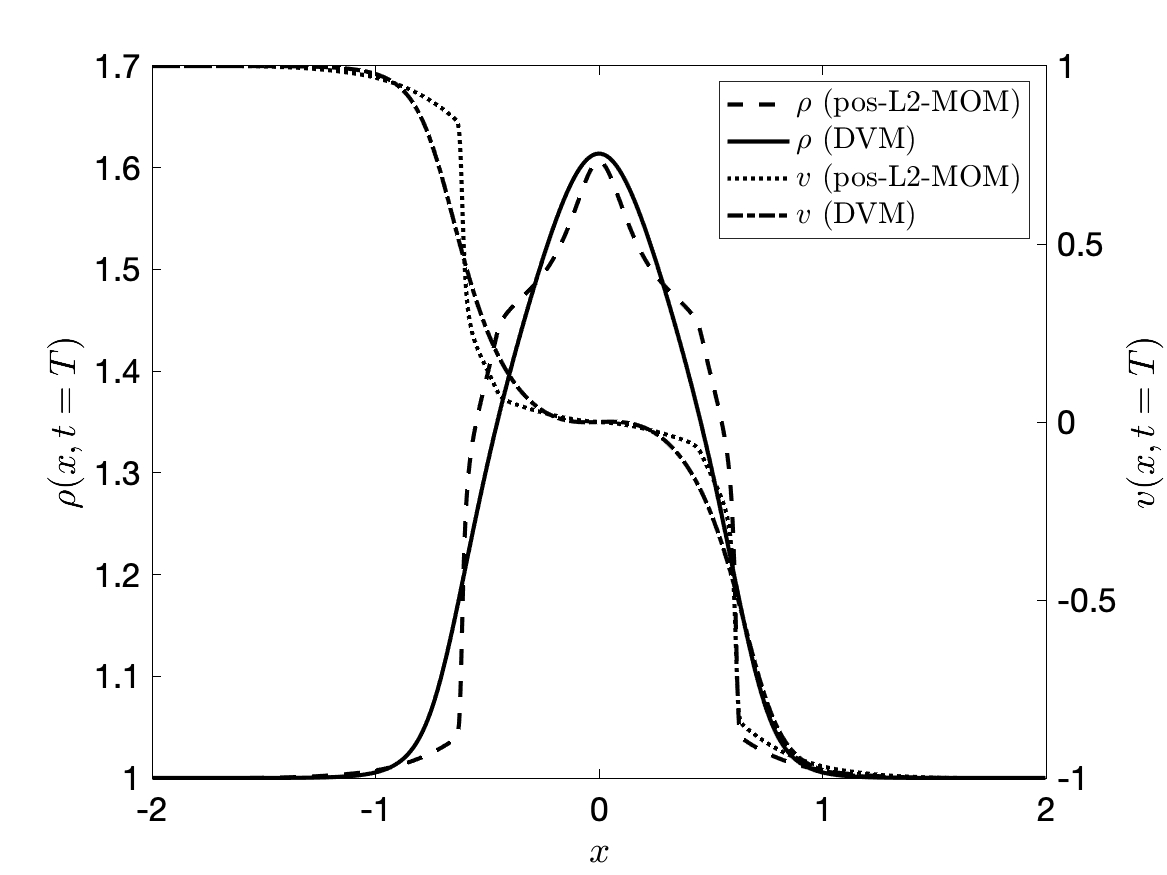}}
\subfloat[$M=7$, $\opn{Kn}=0.1$]{\label{fig: test3 M7_Kn0p1}\includegraphics[width=2.8in]{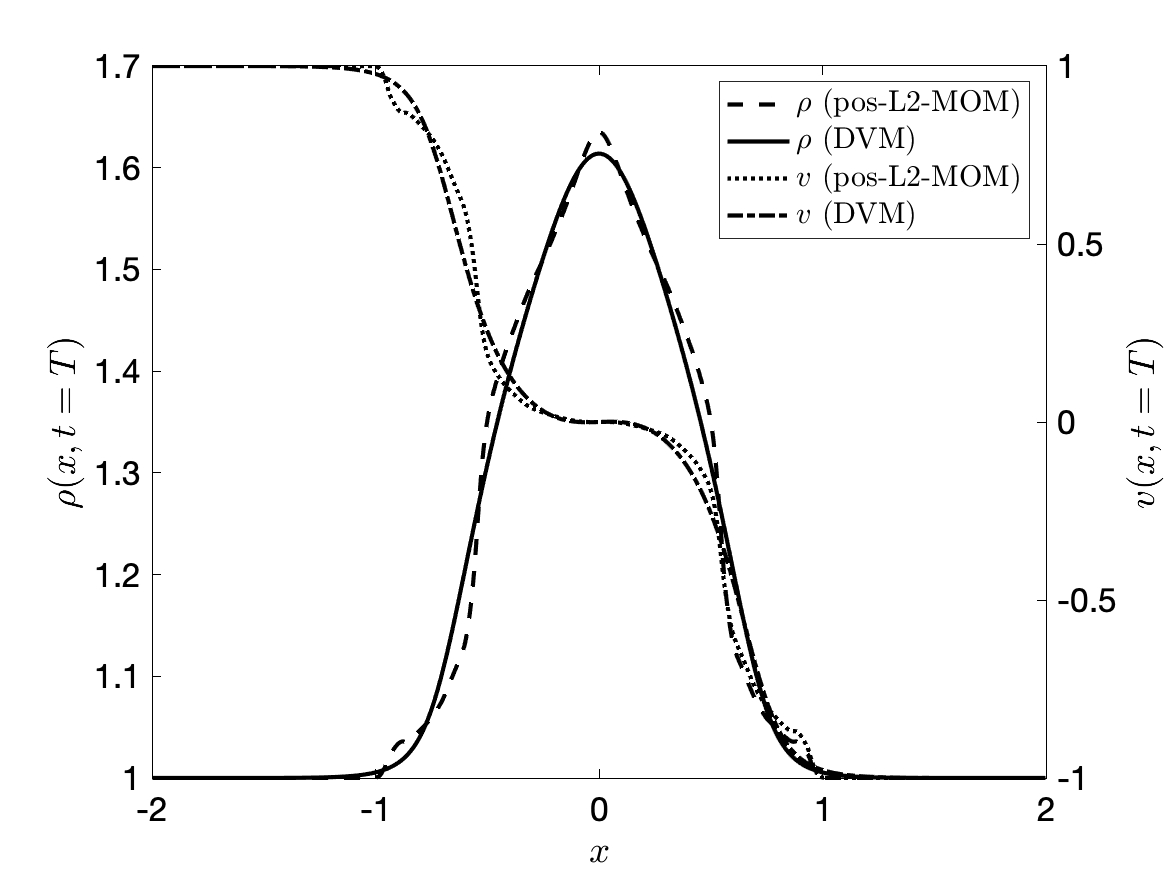}}
\hfill
\subfloat[$M=5$, $\opn{Kn}=0.01$]{\label{fig: test3 M5_Kn0p01}\includegraphics[width=2.8in]{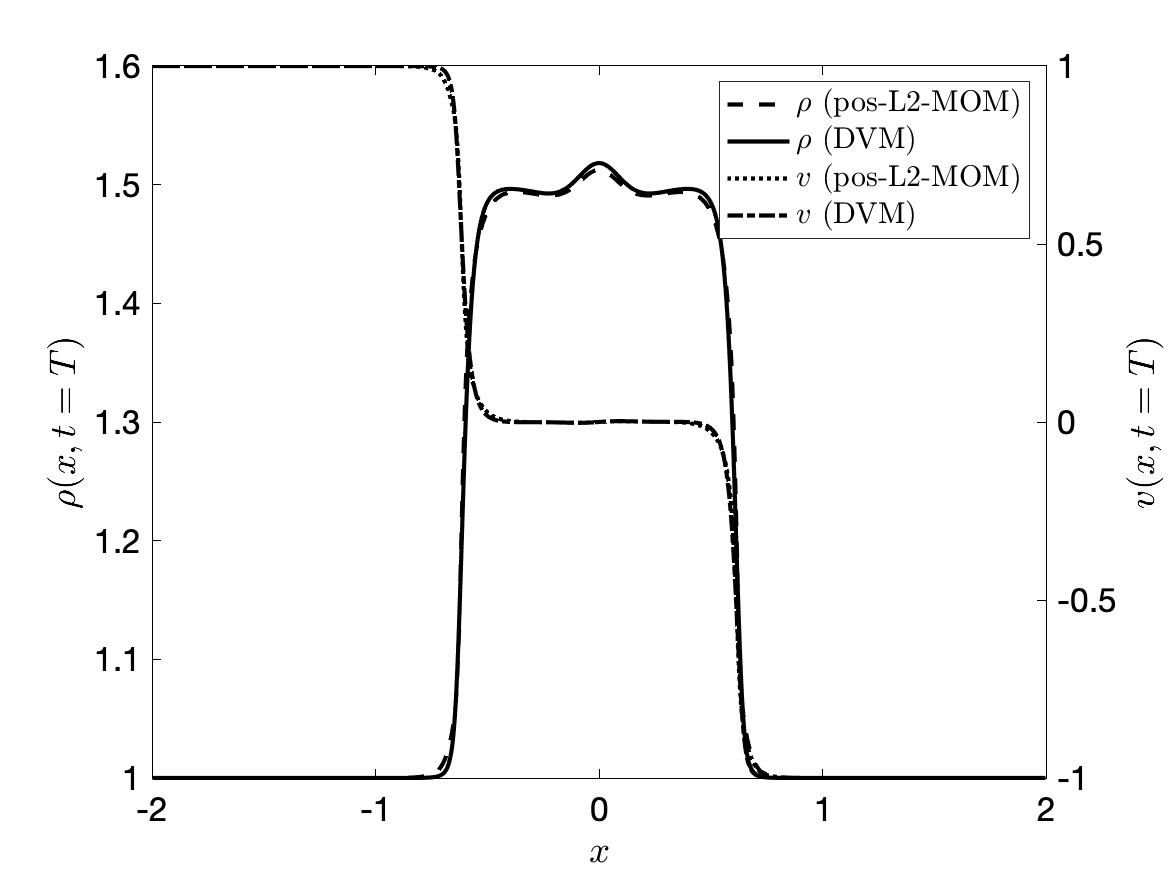}}
\subfloat[$M=7$, $\opn{Kn}=0.01$]{\label{fig: test3 M7_Kn0p01}\includegraphics[width=2.8in]{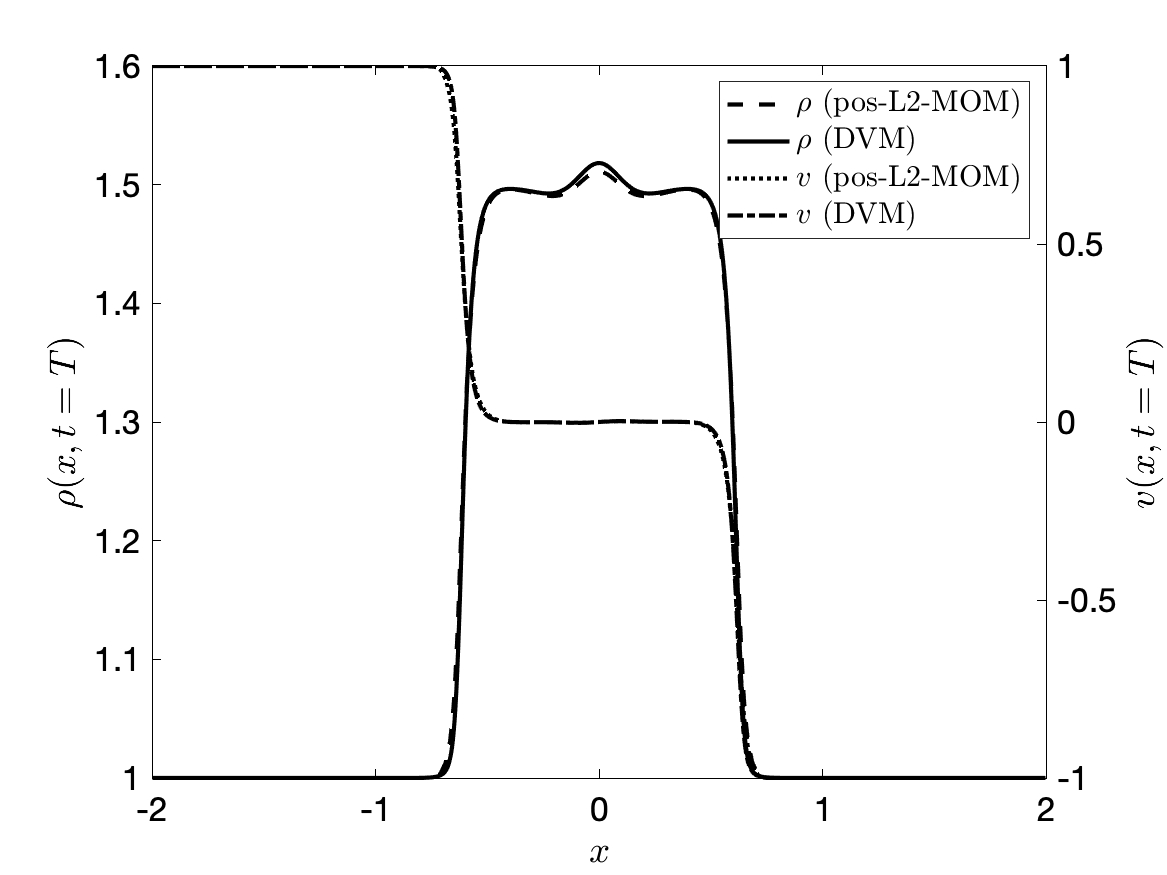}}
\caption{Results for test case-3. Density and velocity profiles for different values of $M$ and different Knudsen numbers. The left and the right y-axis is for density and velocity, respectively. } 

\end{figure}

For $\opn{Kn}=0.1$ and $M=5$, \Cref{fig: test3 M5_Kn0p1} compares the density and the velocity computed using the DVM and the pos-L2-QMOM. The results for temperature are similar and are not shown for brevity. The pos-L2-MOM performs well and results in an error of $\mcal E_{cons}(5,10^3) = 6.8\times 10^{-2}$. Furthermore, increasing the value of $M$ from $5$ to $7$ improves the results and the error reduces to $\mcal E_{cons}(7,10^3) = 2.5\times 10^{-2}$---\Cref{fig: test3 M7_Kn0p1} shows the result for $M=7$. Reducing the Knudsen number to $0.01$, improves the results for both $M=5$ and $M=7$---see \Cref{fig: test3 M5_Kn0p01} and \Cref{fig: test3 M7_Kn0p01}. For both the values of $M$, we obtained an error of $\mcal E_{cons}(5/7,10^3) = 9\times 10^{-3}$, which is approximately $1/3$ of the error for $\opn{Kn} = 0.1$. Note that similar to the previous test case, the error for $\opn{Kn}=0.01$ is dominated by the error in the space-time discretization. Therefore, increasing $M$ from $5$ to $7$ does not offer any improvement.

\subsection{Test case-4}
Under the limited computational resources, we were unable to compute a highly-refined reference solution in multi-dimensions. For this reason, we refrain from performing a convergence study for the present test case. Rather, we compare our moment method to a sufficiently refined DVM and showcase an improvement in the moment solution by increasing $M$. For both the DVM and the moment method, we consider tensorized Gauss-Legendre quadrature points with $N=40$ quadrature points in each direction. We place these quadrature points inside $\Omega_\xi = [\xiMin,\xiMax]\times [\xiMin,\xiMax]$ with $\xiMax = 7$ and $\xiMin = -7$. We discretize the spatial domain with $150\times 150$ uniform elements with grid-size $\Delta x = 1.3\times 10^{-2}$. We consider a constant time-step of $\Delta = \Delta x/(4\times \xiMax)$. 

\begin{figure}[ht!]
\centering
\subfloat[DVM, $\opn{Kn}=0.1$]{\label{fig: test4 rho DVM}\includegraphics[width=2.8in]{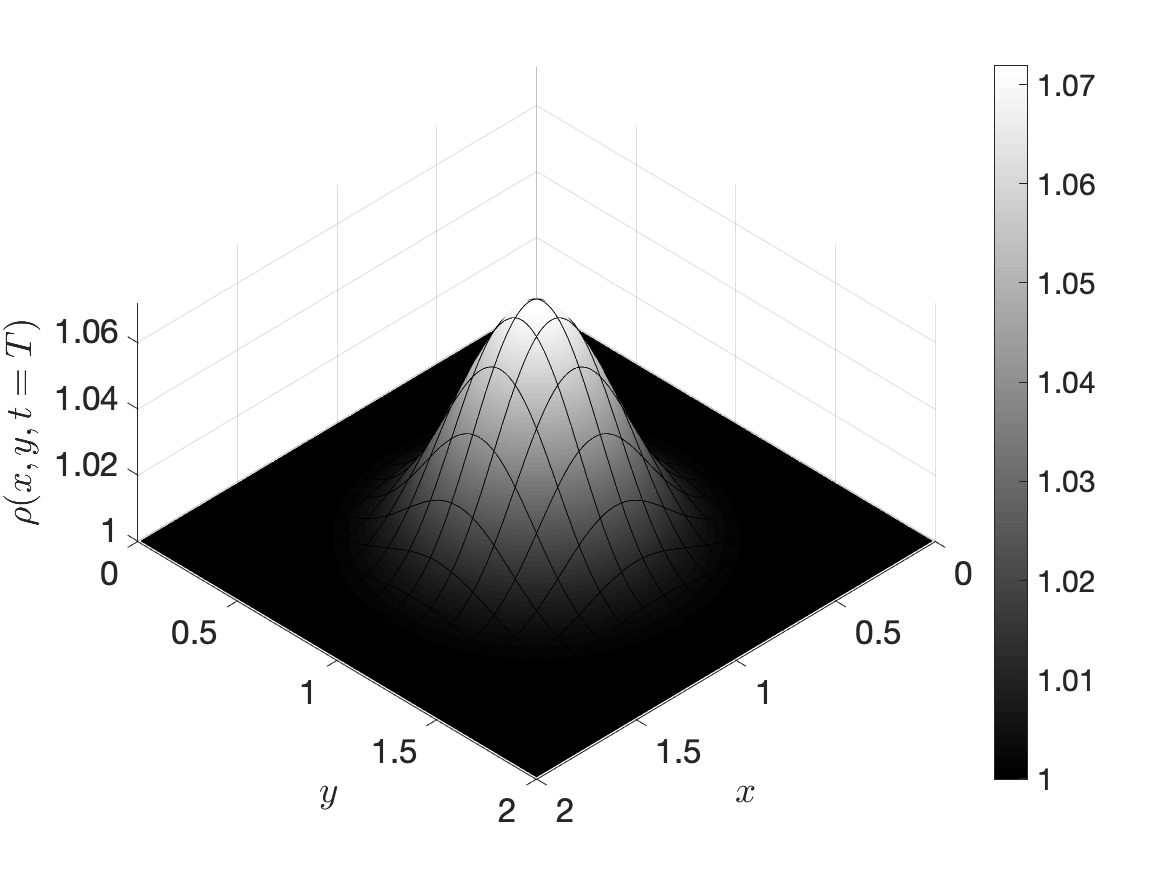}}
\subfloat[$M=3$, $\opn{Kn}=0.1$]{\label{fig: test4 rho M3}\includegraphics[width=2.8in]{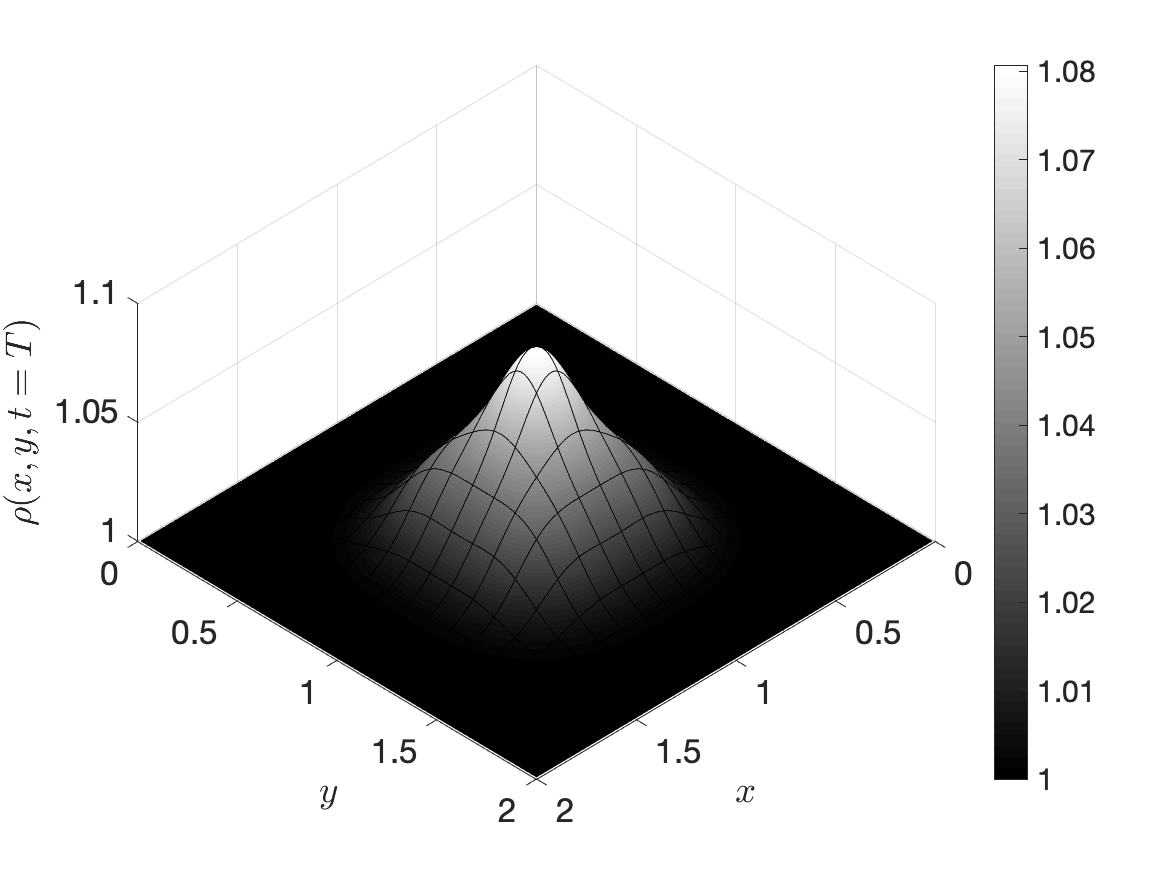}}
\hfill
\subfloat[$M=5$, $\opn{Kn}=0.1$]{\label{fig: test4 rho M5}\includegraphics[width=2.8in]{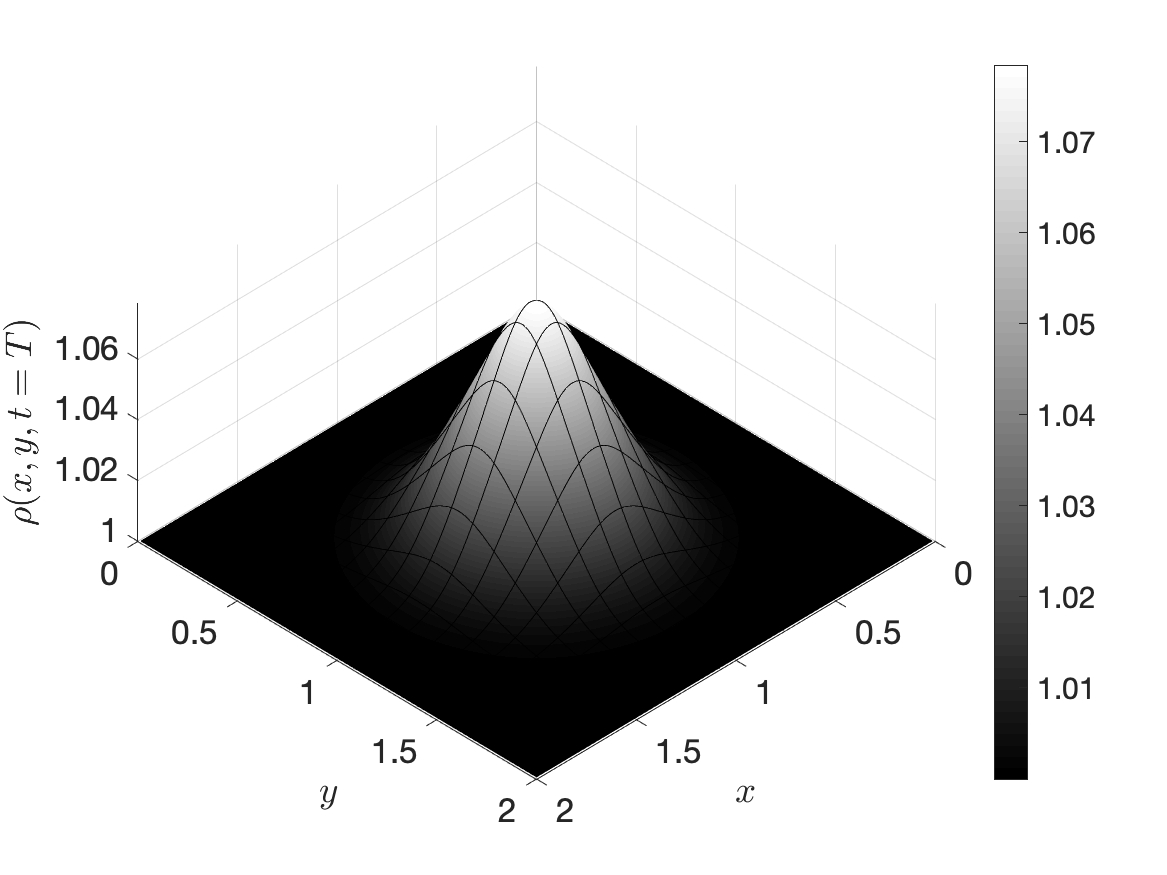}}
\caption{Results for test case-4. Density profiles at $t=T$. } 
\end{figure}

As time progresses, the density disperses into the spatial domain. This is made clear by \Cref{fig: test4 rho DVM} that shows the density profile at $t=T$ computed using the DVM. At the same time-instance, \Cref{fig: test4 rho M3} and \Cref{fig: test4 rho M5} show the density profile at $t=T$ computed using the pos-L2-MOM with $M=3$ and $M=5$, respectively. As expected, both the density profiles are positive. Furthermore, the moment solution appears to improve upon increasing the value of $M$. The improvement is quantified by the decrease in the relative L2-error in density shown in \Cref{table: error}.

\begin{figure}[ht!]
\centering
\subfloat[$v_1(x,y=1,t=T)$, $\opn{Kn}=0.1$]{\label{fig: test4 v1 }\includegraphics[width=2.8in]{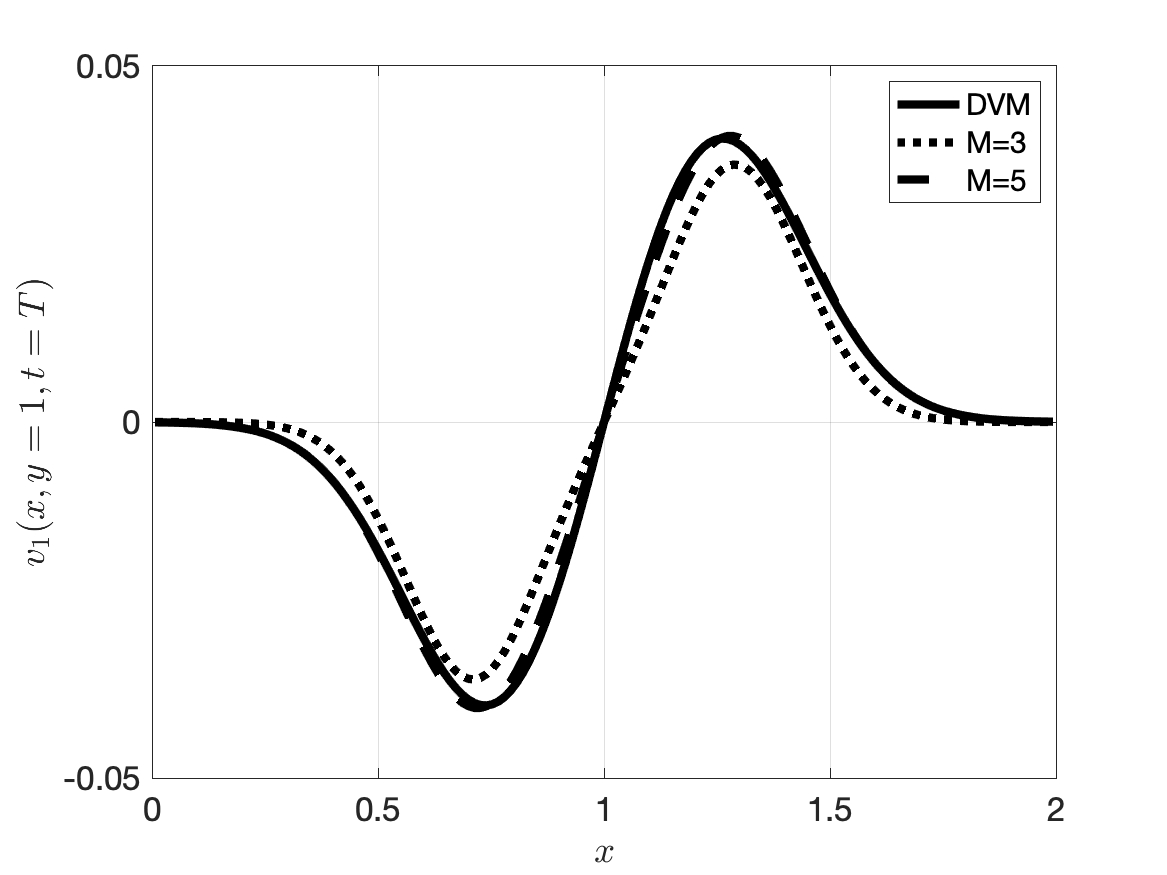}}
\subfloat[$\theta(x,y=1,t=T)$, $\opn{Kn}=0.1$]{\label{fig: test4 theta}\includegraphics[width=2.8in]{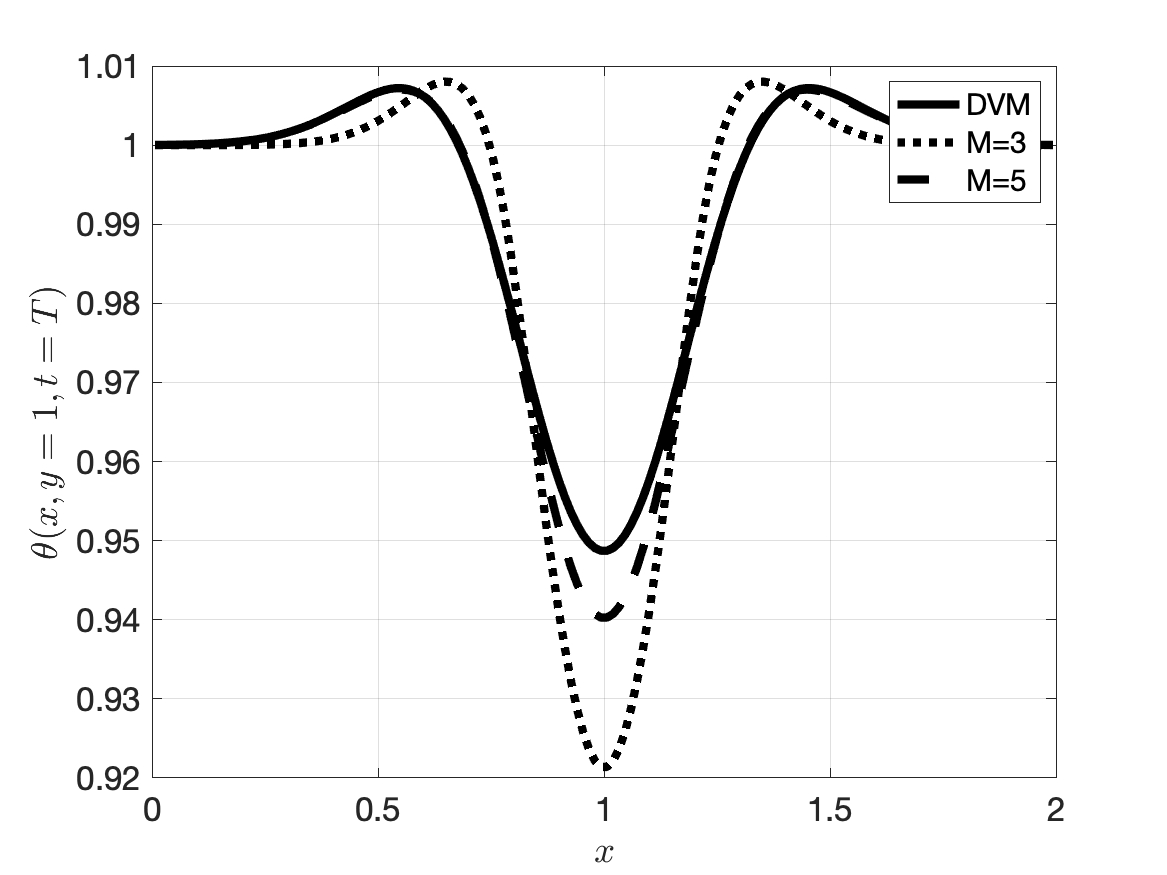}}
\caption{Results for test case-4. $v_1$ and $\theta$ profiles. } \label{fig: test4 v theta}
\end{figure}

The dispersion of the micro-bubble triggers a flow velocity and a temperature gradient. \Cref{fig: test4 v theta} compares the $x_1$ velocity component and the temperature along a cross-section of the spatial domain computed using the pos-L2-MOM and the DVM. The results for the $x_2$ velocity component are similar and are not shown for brevity. As expected, similar to density, the results for both the velocity and the temperature appear to improve as $M$ is increased from $3$ to $5$, the relative L2-error shown in \Cref{table: error} indicates the same. We note that, as compared to the previous test cases, the moment method performs better for the present test case. A possible reason for this could be that our DVM solution is not as refined as for the previous test cases---the previous test cases consider a 1D velocity grid of $350$ points whereas the present test case considers a tensorized grid of $40\times 40$ points.  

\begin{table}[ht!]
\centering
\begin{tabular}{ c | c | c | c | c } 
$M$ & $\rho$ & $v_1$ & $v_2$ & $\theta$ \\ [0.5ex] 
 \hline
 3 & $1.6\times 10^{-3}$ & $2\times 10^{-1}$ & $2.1\times 10^{-1}$ & $2.8\times 10^{-3}$ \\ 
 \hline
 5 & $5.3\times 10^{-4}$ & $5\times 10^{-2}$ & $4.8\times 10^{-2}$ & $5.4\times 10^{-4}$ 
\end{tabular}
\caption{Results for test case-4. Relative $L^2(\Omega)$-error in different macroscopic quantities at $t=T$ and $Kn=0.1$. }
\label{table: error}
\end{table}
\section{Conclusions}
We proposed a positive moment method for the Boltzmann-BGK equation based upon L2-minimization. We showed that on a space-time discrete level both the feasibility of the minimization problem and the stability of the moment approximation can be ensured via a CFL-type condition. Our proof of booth these properties relied on relating our moment method to a discrete-velocity-method. Through a entropy-minimization based discretization of the collision operator, we ensured that our moment approximation conserves mass, momentum and energy. We also extended our method to a multi-dimensional space-velocity domain. With the help of numerical experiments, we studied the accuracy of our method for both single and multi-dimensional space-velocity domains. Our method performed well for a broad range of problems involving strong shocks, beam interaction and micro-bubble dispersion, and retained accuracy for a broad range of Knudsen numbers.

\bibliographystyle{abbrv}
\bibliography{papers}
\end{document}